\documentclass[12pt,a4paper]{article}
\usepackage{}
\usepackage{amsfonts}
\usepackage{amssymb}
\usepackage{mathrsfs}
\usepackage{amsmath}
\usepackage{amsthm}
\usepackage{enumerate}
\usepackage{cite}
\usepackage[all]{xy}
\allowdisplaybreaks
\newtheorem{defn}{Definition}[section]
\newtheorem{thm}[defn]{Theorem}
\newtheorem{lem}[defn]{Lemma}
\newtheorem{prop}[defn]{Proposition}
\newtheorem{cor}[defn]{Corollary}
\newtheorem{eg}[defn]{Example}
\newtheorem{re}[defn]{Remark}

\newcommand\relphantom[1]{\mathrel{\phantom{#1}}}

\newcommand{\bdefn}{\begin{defn}}
\newcommand{\edefn}{\end{defn}}
\newcommand{\bthm}{\begin{thm}}
\newcommand{\ethm}{\end{thm}}
\newcommand{\blem}{\begin{lem}}
\newcommand{\elem}{\end{lem}}
\newcommand{\bprop}{\begin{prop}}
\newcommand{\eprop}{\end{prop}}
\newcommand{\bcor}{\begin{cor}}
\newcommand{\ecor}{\end{cor}}
\newcommand{\beg}{\begin{eg}}
\newcommand{\eeg}{\end{eg}}
\newcommand{\bre}{\begin{re}}
\newcommand{\ere}{\end{re}}
\newcommand{\bpf}{\begin{proof}}
\newcommand{\epf}{\end{proof}}

\newcommand{\ch}{{\rm ch}}
\newcommand{\id}{{\rm id}}
\newcommand{\ad}{{\rm ad}}
\newcommand{\End}{{\rm End}}
\newcommand{\K}{\mathbb{K}}

\newcommand{\Z}{\mathbb{Z}}
\newcommand{\g}{\mathfrak{g}}
\newcommand{\A}{\mathscr{A}}
\newcommand{\B}{\mathscr{B}}
\newcommand{\x}{\mathscr{X}}
\newcommand{\y}{\mathscr{Y}}
\newcommand{\z}{\mathscr{Z}}
\newcommand{\F}{\mathscr{F}}
\newcommand{\G}{\mathscr{G}}
\newcommand{\XY}{\mathscr{X}\cdot\mathscr{Y}}
\newcommand{\yx}{\mathscr{Y}\cdot\mathscr{X}}
\newcommand{\xz}{\mathscr{X}\cdot\mathscr{Z}}

\newcommand{\yz}{\mathscr{Y}\cdot\mathscr{Z}}
\newcommand{\xij}{\mathscr{X}_i\cdot\mathscr{X}_j}
\newcommand{\xst}{\mathscr{X}_s\cdot\mathscr{X}_t}

\newcommand{\benu}{\begin{enumerate}}
\newcommand{\eenu}{\end{enumerate}}
\newcommand{\bc}{\begin{center}}
\newcommand{\ec}{\end{center}}
\newcommand{\bea}{\begin{eqnarray}}
\newcommand{\eea}{\end{eqnarray}}
\newcommand{\Bea}{\begin{eqnarray*}}
\newcommand{\Eea}{\end{eqnarray*}}
\newcommand{\beq}{\begin{equation}}
\newcommand{\eeq}{\end{equation}}
\newcommand{\Beq}{\begin{equation*}}
\newcommand{\Eeq}{\end{equation*}}
\newcommand{\bspl}{\begin{split}}
\newcommand{\espl}{\end{split}}

\newcommand{\supercite}[1]{\textsuperscript{\cite{#1}}}

\begin{document}
\title{{\bf  On the cohomology and extensions of first-class $n$-Lie superalgebras}}
\author{\normalsize \bf Yao Ma,  Liangyun Chen}
\date{{{\small{ School of Mathematics and Statistics,  Northeast Normal University, Changchun 130024, China
 }}}} \maketitle
\date{}

\begin{abstract}

An $n$-Lie superalgebra of parity 0 is called a first-class $n$-Lie superalgebra. In this paper, we give the representation and cohomology for a first-class $n$-Lie superalgebra and obtain a relation between extensions of a first-class $n$-Lie superalgebra $\mathfrak{b}$ by an abelian one $\mathfrak{a}$ and $Z^1(\mathfrak{b}, \mathfrak{a})_{\bar{0}}$. We also introduce the notion of $T^*$-extensions of first-class $n$-Lie superalgebras and prove that every finite-dimensional nilpotent metric first-class $n$-Lie superalgebra $(\g,\langle ,\rangle_{\g})$ over  an algebraically closed field of characteristic not 2 is isometric to a suitable $T^*$-extension.

\bigskip

\end{abstract}
\renewcommand{\thefootnote}{\fnsymbol{footnote}}
\footnote[0]{ Corresponding author(L. Chen): chenly640@nenu.edu.cn.}
\footnote[0]{Supported by  NNSF of China (No.11171055),  NSF of  Jilin province (No.201115006), Scientific
Research Foundation for Returned Scholars
    Ministry of Education of China. }

\section{Introduction}

V. T. Filippov introduced the definition of $n$-Lie algebras in 1985, and a structure theory of finite-dimensional $n$-Lie algebras over a field $\K$ of characteristic 0 was developed\supercite{F,K2,L}. $n$-Lie algebras were found useful in the research for M2-branes in the string theory and were closely linked to the Pl\"{u}cker relation in literature in physics\supercite{G,BL1,BL2,P}. $n$-Lie superalgebras are more general structures containing $n$-Lie algebras and Lie superalgebras, whose definition was introduced by N. Cantarini and V.G. Kac\supercite{CK}.
Cohomologies are powerful tools in mathematics, which can be applied to algebras and topologies as well as the theory of smooth manifolds or of holomorphic functions. The cohomology of Lie algebras was defined by C. Chevalley and S. Eilenberg in order to give an algebraic construction of the cohomology of the underlying topological spaces of compact Lie groups\supercite{CE}. The cohomology of Lie superalgebras was introduced by M. Scheunert and R. B. Zhang\supercite{SZ} and was used in mathematics and theoretical physics\supercite{BMPZ,K3}: the theory of cobordisms, invariant differential operators, central extensions and deformations, etc. The theory of cohomology for $n$-Lie algebras and Lie triple systems can be found in \cite{AI,Y}. This paper discusses first-class $n$-Lie superalgebras, i.e., $n$-Lie superalgebras of parity 0, and gives the cohomology of first-class $n$-Lie superalgebras.

The extension is an important way to find a larger algebra and there are many extensions such as double extensions and Kac-Moody extensions, etc. In 1997, Bordemann introduced the notion of $T^*$-extensions of Lie algebras\supercite{B} and proved that every nilpotent finite-dimensional algebra over an algebraically closed field carrying a nondegenerate invariant symmetric bilinear form is a suitable $T^*$-extension. The method of $T^*$-extension was used in \cite{BBM,MR,K1} and was generalized to many other algebras recently\supercite{LWD,BL3,LZ,LZ2}. This paper researches general extensions and $T^*$-extensions of first-class $n$-Lie superalgebras.

This paper is organized as follows. In section 2, we give the representation and the cohomology for a first-class $n$-Lie superalgebra. In section 3, we give a one-to-one correspondence between extensions of a first-class $n$-Lie superalgebra $\mathfrak{b}$ by an abelian one $\mathfrak{a}$ and $Z^1(\mathfrak{b}, \mathfrak{a})_{\bar{0}}$. In section 4, we introduce the notion of $T^*$-extensions of first-class $n$-Lie superalgebras and prove that every finite-dimensional nilpotent metric first-class $n$-Lie superalgebra $(\g,\langle ,\rangle_{\g})$ over  an algebraically closed field of characteristic not 2 is isometric to (a nondegenerate ideal of codimension 1 of) a $T^*$-extension of a nilpotent first-class $n$-Lie superalgebra whose nilpotent length is at most a half of the nilpotent length of $\g$.

\bdefn\supercite{CK}
A $\Z_2$-graded vector space $\g=\g_{\bar{0}}\oplus \g_{\bar{1}}$ is called an $n$-Lie superalgebra of parity $\alpha$, if there is a multilinear mapping $[ ,\cdots, ]: \underbrace{\g\times \cdots \times \g}_n\rightarrow \g$ such that
\begin{align}|[x_1,\cdots, x_n]|=&|x_1|+\cdots+|x_n|+\alpha,\label{eq:1}\\
[x_1,\cdots, x_i, x_{i+1},\cdots, x_n]=&-(-1)^{|x_i| |x_{i+1}|}[x_1,\cdots, x_{i+1}, x_i,\cdots, x_n],\label{eq:2}\\
\begin{split}
[x_1,\cdots, x_{n-1}, [y_1,\cdots, y_n]]=&(-1)^{\alpha(|x_1|+\cdots+|x_{n-1}|)}\sum_{i=1}^n(-1)^{(|x_1|+\cdots+|x_{n-1}|)(|y_1|+\cdots+|y_{i-1}|)}\label{eq:3}\\
&\relphantom{=}\cdot[y_1,\cdots,[x_1,\cdots, x_{n-1},y_i],\cdots,y_n],
\end{split}\end{align}
where $\alpha\in \Z_2$ and $|x|\in\Z_2$ denotes the degree of a homogeneous element $x\in\g$.
%\begin{align}
%&|[a_1,\cdots, a_n]|=|a_1|+\cdots+|a_n|+\alpha,\\
%&[a_1,\cdots, a_i, a_{i+1},\cdots, a_n]=-(-1)^{|a_i| |a_{i+1}|}[a_1,\cdots, a_{i+1}, a_i,\cdots, a_n],\\
%\begin{split}
%&[a_1,\cdots, a_{n-1}, [b_1,\cdots, b_n]]\\
%&=(-1)^{\alpha(|a_1|+\cdots+|a_{n-1}|)}\sum_{i=1}^n(-1)^{(|a_1|+\cdots+|a_{n-1}|)(|b_1|+\cdots+|b_{i-1}|)}[b_1,\cdots,[a_1,\cdots, a_{n-1},b_i],\cdots,b_n],
%\end{split}
%\end{align}
\edefn

In this paper, we only consider the case ``$\alpha=0$'', then Equations (\ref{eq:1}) and (\ref{eq:3}) can be read:
\begin{align}|[x_1,\cdots, x_n]|=&|x_1|+\cdots+|x_n|,\tag{\ref{eq:1}$'$}\label{eq:1'}\\
\begin{split}\label{eq:3'}
[x_1,\cdots, x_{n-1}, [y_1,\cdots, y_n]]=&\sum_{i=1}^n(-1)^{(|x_1|+\cdots+|x_{n-1}|)(|y_1|+\cdots+|y_{i-1}|)}\\
&\relphantom{=}\cdot[y_1,\cdots,[x_1,\cdots, x_{n-1},y_i],\cdots,y_n].
\end{split}\tag{\ref{eq:3}$'$}\end{align}
%\beq  \label{eq:1'} |[a_1,\cdots, a_n]|=|a_1|+\cdots+|a_n|, \tag{\ref{eq:1}$'$}\eeq
%\beq\label{eq:3'}\begin{split}
%&[a_1,\cdots, a_{n-1}, [b_1,\cdots, b_n]]\\
%=&\sum_{i=1}^n(-1)^{(|a_1|+\cdots+|a_{n-1}|)(|b_1|+\cdots+|b_{i-1}|)}[b_1,\cdots,[a_1,\cdots, a_{n-1},b_i],\cdots,b_n].
%\end{split}\tag{\ref{eq:3}$'$}\eeq
We call an $n$-Lie superalgebra of parity 0 \textit{a first-class $n$-Lie superalgebra}. It is clear that $n$-Lie algebras and Lie superalgebras are contained in first-class $n$-Lie superalgebras. In the sequel, when the notation ``$|x|$'' appears, it means that $x$ is a homogeneous element of degree $|x|$.

\section{Cohomology for first-class $n$-Lie superalgebras}

\bdefn
Let $\g$ be a first-class $n$-Lie superalgebra. $\x=x_1\wedge\cdots\wedge x_{n-1}\in \g^{\wedge^{n-1}}$ is called a fundamental object of $\g$ and $\forall z\in\g, \x\cdot z:=[x_1,\cdots, x_{n-1},z]$. Then a fundamental object defines an inner derivation of $\g$ and $|\x|=|x_1|+\cdots+|x_{n-1}|$.
\edefn

Let $\x=x_1\wedge\cdots\wedge x_{n-1}$ and $\y=y_1\wedge\cdots\wedge y_{n-1}$ be two fundamental objects of $\g$. The  composition $\XY\in \g^{\wedge^{n-1}}$ is defined by
$$\XY=\sum_{i=1}^{n-1}(-1)^{|\x|(|y_1|+\cdots+|y_{i-1}|)}y_1\wedge\cdots\wedge\x\cdot y_i\wedge\cdots\wedge y_{n-1}.$$
Then $(\XY)\cdot z=\sum\limits_{i=1}^{n-1}(-1)^{|\x|(|y_1|+\cdots+|y_{i-1}|)}[y_1, \cdots,\x\cdot y_i,\cdots, y_{n-1},z]$.

\bprop\label{prop:1}
Suppose that $\x=x_1\wedge\cdots\wedge x_{n-1}$, $\y=y_1\wedge\cdots\wedge y_{n-1}$ and $\z=z_1\wedge\cdots\wedge z_{n-1}$ are fundamental objects of $\g$ and $z$ is an arbitrary element in $\g$. Then \begin{align}
\x\cdot(\y\cdot z)&=(\XY)\cdot z+(-1)^{|\x||\y|}\y\cdot(\x\cdot z),\label{eq:4}\\
\x\cdot(\yz)&=(\XY)\cdot\z+(-1)^{|\x||\y|}\y\cdot(\xz),\label{eq:5}\\
(\XY)\cdot z&=-(-1)^{|\x||\y|}(\yx)\cdot z.\label{eq:6}
\end{align}
\eprop
\bpf It's easy to see that (\ref{eq:4}) is equivalent to (\ref{eq:3'}). Note that
\begin{align}
&\x\cdot(\yz)=\x\cdot\left(\sum_{i=1}^{n-1}(-1)^{|\y|(|z_1|+\cdots+|z_{i-1}|)}z_1\wedge\cdots\wedge\y\cdot z_i\wedge\cdots\wedge z_{n-1}\right)\notag\\
=&\sum_{i=1\atop i< j}^{n-1}(-1)^{|\y|(|z_1|+\cdots+|z_{i-1}|)}(-1)^{|\x|(|z_1|+\cdots+|z_{j-1}|+|\y|)}z_1\wedge\!\cdots\!\wedge\y\cdot z_i\wedge\!\cdots\!\wedge\x\cdot z_j\wedge\!\cdots\!\wedge z_{n-1}\tag{\ref{eq:5}a}\label{eq:5a}\\
&+\sum_{i=1\atop j< i}^{n-1}(-1)^{|\y|(|z_1|+\cdots+|z_{i-1}|)}(-1)^{|\x|(|z_1|+\cdots+|z_{j-1}|)}z_1\wedge\!\cdots\!\wedge\x\cdot z_j\wedge\!\cdots\!\wedge\y\cdot z_i\wedge\!\cdots\!\wedge z_{n-1}\tag{\ref{eq:5}b}\label{eq:5b}\\
&+\sum_{i=1}^{n-1}(-1)^{(|\x|+|\y|)(|z_1|+\cdots+|z_{i-1}|)}z_1\wedge\cdots\wedge\x\cdot(\y\cdot z_i)\wedge\cdots\wedge z_{n-1}\tag{\ref{eq:5}c}\label{eq:5c}.
\end{align}
Similarly,
\begin{align}
&\y\cdot(\xz)\notag\\
=&\sum_{i=1\atop i<j}^{n-1}(-1)^{|\x|(|z_1|+\cdots+|z_{i-1}|)}(-1)^{|\y|(|z_1|+\cdots+|z_{j-1}|+|\x|)}z_1\wedge\!\cdots\!\wedge\x\cdot z_i\wedge\!\cdots\!\wedge\y\cdot z_j\wedge\!\cdots\!\wedge z_{n-1}\tag{\ref{eq:5}d}\label{eq:5d}\\
&+\sum_{i=1\atop j< i}^{n-1}(-1)^{|\x|(|z_1|+\cdots+|z_{i-1}|)}(-1)^{|\y|(|z_1|+\cdots+|z_{j-1}|)}z_1\wedge\cdots\wedge\y\cdot z_j\wedge\!\cdots\!\wedge\x\cdot z_i\wedge\!\cdots\!\wedge z_{n-1}\tag{\ref{eq:5}e}\label{eq:5e}\\
&+\sum_{i=1}^{n-1}(-1)^{(|\y|+|\x|)(|z_1|+\cdots+|z_{i-1}|)}z_1\wedge\cdots\wedge\y\cdot(\x\cdot z_i)\wedge\cdots\wedge z_{n-1}\tag{\ref{eq:5}f}\label{eq:5f},
\end{align}
and
\beq
(\XY)\cdot\z=\sum_{i=1}^{n-1}(-1)^{(|\x|+|\y|)(|z_1|+\cdots+|z_{i-1}|)}z_1\wedge\cdots\wedge(\XY)\cdot z_i\wedge\cdots\wedge z_{n-1}\tag{\ref{eq:5}g}\label{eq:5g}.
\eeq

It can be checked that (\ref{eq:5a})+(\ref{eq:5b})=$(-1)^{|\x||\y|}$(\ref{eq:5d})+$(-1)^{|\x||\y|}$(\ref{eq:5e}). By (\ref{eq:4}), we conclude (\ref{eq:5c})=(\ref{eq:5g})+$(-1)^{|\x||\y|}$(\ref{eq:5f}). Thus (\ref{eq:5}) holds.

Use (\ref{eq:4}), by exchanging $\x$ and $\y$, we have
\beq\y\cdot(\x\cdot z)=(\yx)\cdot z+(-1)^{|\x||\y|}\x\cdot(\y\cdot z).\label{eq:7}\eeq
Compare (\ref{eq:4}) with (\ref{eq:7}) we obtain (\ref{eq:6}).
\epf

\bdefn
Let $\g$ be a first-class $n$-Lie superalgebra and $V=V_{\bar{0}}\oplus V_{\bar{1}}$ be a $\Z_2$-graded vector space over a field $\K$. A graded representation $\rho$ of $\g$ on $V$ is a linear map $\rho: \g^{\wedge^{n-1}}\rightarrow \End(V), \x\mapsto\rho(\x)=\rho(x_1,\cdots,x_{n-1})$ such that
\begin{align}
\rho(\x)\cdot V_\alpha\subseteq& V_{\alpha+|\x|}, \forall \alpha\in \Z_2,\\
\rho(\x)\rho(\y)=&\rho(\XY)+(-1)^{|\x||\y|}\rho(\y)\rho(\x),\\
\begin{split}\label{eq:rho3}
\rho(x_1,\cdots,x_{n-2},[y_1,\cdots,y_n])=
&\sum_{i=1}^n(-1)^{n-i}(-1)^{(|x_1|+\cdots+|x_{n-2}|)(|y_1|+\cdots+\widehat{|y_i|}+\cdots+|y_n|)}\\
%&\relphantom{=}\cdot(-1)^{|y_i|(|y_{i+1}+\cdots+|y_n|)}\\
&\relphantom{=}\cdot(-1)^{|y_i|(|y_{i+1}|+\cdots+|y_n|)}\rho(y_1,\!\cdots\!,\widehat{y_i},\!\cdots\!,y_n)\rho(x_1,\!\cdots\!,x_{n-2},y_i),\end{split}
\end{align}
for all $\x,\y\in\g^{\wedge^{n-1}}$ and $x_1,\cdots,x_{n-2},y_1,\cdots,y_n\in\g$, and the sign \textasciicircum ~ indicates that the element below it must be omitted.
The $\Z_2$-graded representation space $V$ is said to be a graded $\g$-module.
\edefn

If we use a supersymmetric notation $[x_1,\cdots,x_{n-1}, v]$(like (\ref{eq:2})) to denote $\rho(\x)\cdot v$ and set $[x_1,\cdots,x_{n-2}, v_1,v_2]=0$, then $\g\oplus V$ becomes a first-class $n$-Lie superalgebra such that $V$ is a $\Z_2$-graded abelian ideal of $\g$, that is,
$$[V,\underbrace{\g,\cdots,\g}_{n-1}]\subseteq V \quad \text{and}\quad [V,V,\underbrace{\g,\cdots,\g}_{n-2}]=0.$$
In the sequel, we usually abbreviate $\rho(\x)\cdot v$ to $\x\cdot v$.

\beg
Let $\g^*$ be the dual $\Z_2$-graded vector space of a first-class $n$-Lie superalgebra $\g$. Then $\g^*$ is a graded $\g$-module with the coadjoint graded representation $\ad^*:\g^{\wedge^{n-1}}\rightarrow \End(\g^*)$ defined by
$$\ad^*(\x)(f)(z)=-(-1)^{|\x||f|}f(\x\cdot z),$$
for all $\x\in \g^{\wedge^{n-1}}, f\in\g^*$ and $z\in\g$. Moreover, the field $\K$ is a trivial graded $\g$-module.
\eeg

\bdefn
Suppose that $V=V_{\bar{0}}\oplus V_{\bar{1}}$ is a graded $\g$-module. Let
$$C^m(\g, V)={\rm Hom} (\underbrace{\g^{\wedge^{n-1}}\otimes\cdots\otimes\g^{\wedge^{n-1}}}_m\wedge\g, V)$$
denote the set of all $m$-supercochains, $\forall m\geq0$. Then $C^m(\g, V)$ is a $\Z_2$-graded vector space with $C^m(\g, V)_\alpha=\{f\in C^m(\g, V)||f|=\alpha\in\Z_2\}$.
\edefn
%\bprop
%$C^m(\g, V)$ is a graded $\g$-module.
%\eprop
%\bpf
%Let $\rho:\g^{\wedge^{n-1}}\rightarrow \End(C^m(\g, V))$ be defined as
%\Beq\begin{split}
%&(\rho(\y)\cdot f)(\x_1,\cdots,\x_m,z)\\
%=&\y\cdot f(\x_1,\cdots,\x_m,z)-(-1)^{|\y|(|f|+|\x_1|+\cdots+|\x_{m}|)}f(\x_1,\cdots,\x_m,\y\cdot z)\\
%&-\sum_{i=1}^m(-1)^{|\y|(|f|+|\x_1|+\cdots+|\x_{i-1}|)}f(\x_1,\cdots,\yx_i,\cdots,\x_m,z),
%\end{split}\Eeq
%where $\x_i=\x_i^1\wedge\cdots\wedge\x_i^{n-1}\in\g^{\wedge^{n-1}}$, for all $1\leq i\leq m$. Then $\rho$ can be proved to be a graded representation of $\g$ on $C^m(\g, V)$ by Proposition \ref{prop:1} and the fact that
%\beq\begin{split}
%[x_1,\cdots,x_{n-2},[y_1,\cdots,y_n],z]=&\sum_{i=1}^n(-1)^{n-i}(-1)^{(|x_1|+\cdots+|x_{n-2}|)(|y_1|+\cdots+\widehat{|y_i|}+\cdots+|y_n|)}\\
%&\relphantom{=}\cdot[y_1,\cdots,\widehat{y_i},\cdots,y_n,[x_1,\cdots,x_{n-2},y_i,z]],
%\end{split}\eeq
%for all $x_i$, $y_j$ and $z\in\g$ ($1\leq i\leq n-2, 1\leq j\leq n$).
%\epf

\bdefn\label{def:delta}
We now define a linear map $\delta: C^m(\g, V)\rightarrow C^{m+1}(\g, V)$ by
\begin{align*}
&(\delta f)(\x_1,\cdots, \x_{m+1}, z)\\
=&\sum_{i<j}(-1)^i(-1)^{|\x_i|(|\x_{i+1}|+\cdots+|\x_{j-1}|)}f(\x_1,\cdots,\widehat{\x_i},\cdots,\xij,\cdots,\x_{m+1},z)\\
&+\sum_{i=1}^{m+1}(-1)^i(-1)^{|\x_i|(|\x_{i+1}|+\cdots+|\x_{m+1}|)}f(\x_1,\cdots,\widehat{\x_i},\cdots,\x_{m+1},\x_i\cdot z)\\
& +\sum_{i=1}^{m+1}(-1)^{i+1}(-1)^{|\x_i|(|f|+|\x_{1}|+\cdots+|\x_{i-1}|)}\x_i\cdot f(\x_1,\cdots,\widehat{\x_i},\cdots,\x_{m+1}, z)\\
&  +(-1)^m(f(\x_1,\cdots,\x_m, ~~)\cdot \x_{m+1})\cdot z,
\end{align*}
where $\x_i=\x_i^1\wedge\cdots\wedge\x_i^{n-1}\in\g^{\wedge^{n-1}}, i=1,\cdots,m+1, z\in\g$ and the last term is defined by
\Beq\begin{split}
(f(\x_1,\cdots,\x_m,~~ )\cdot \x_{m+1})\cdot z=&\sum_{i=1}^{n-1}(-1)^{(|f|+|\x_{1}|+\cdots+|\x_{m}|)(|\x_{m+1}^1|+\cdots+|\x_{m+1}^{i-1}|)}\\
&\relphantom{=}\cdot[\x_{m+1}^1,\cdots,f(\x_1,\cdots,\x_m,\x_{m+1}^i),\cdots,\x_{m+1}^{n-1},z].
\end{split}\Eeq
\edefn

We now check that $\delta^2=0$. In fact,
\begin{align}
&(\delta^2 f)(\x_1,\cdots, \x_{m+2}, z)\notag\\
=&\sum_{i<j}(-1)^i(-1)^{|\x_i|(|\x_{i+1}|+\cdots+|\x_{j-1}|)}\delta f(\x_1,\cdots,\widehat{\x_i},\cdots,\xij,\cdots,\x_{m+2},z)\notag\\
&+\sum_{i=1}^{m+2}(-1)^i(-1)^{|\x_i|(|\x_{i+1}|+\cdots+|\x_{m+2}|)}\delta f(\x_1,\cdots,\widehat{\x_i},\cdots,\x_{m+2},\x_i\cdot z)\notag\\
& +\sum_{i=1}^{m+2}(-1)^{i+1}(-1)^{|\x_i|(|f|+|\x_{1}|+\cdots+|\x_{i-1}|)}\x_i\cdot\delta f(\x_1,\cdots,\widehat{\x_i},\cdots,\x_{m+2}, z)\notag\\
&  +(-1)^{m+1}(\delta f(\x_1,\cdots,\x_{m+1}, ~~)\cdot \x_{m+2})\cdot z,\notag\\
=&\sum_{s<t<i<j}a_{ijst}f(\x_1,\cdots,\widehat{\x_s},\cdots,\xst,\cdots, \widehat{\x_i}, \cdots, \xij,\cdots, \x_{m+2},z)\tag{a1}\\
  &+\sum_{s<i<t<j}\widetilde{a_{ijst}}f(\x_1,\cdots,\widehat{\x_s},\cdots, \widehat{\x_i},\cdots,\xst, \cdots, \xij,\cdots, \x_{m+2},z)\tag{a2}\\
   &+\sum_{s<i<j<t}a_{ijst}f(\x_1,\cdots,\widehat{\x_s},\cdots, \widehat{\x_i},\cdots,\xij, \cdots, \xst,\cdots, \x_{m+2},z)\tag{a3}\\
    &-\sum_{i<s<t<j}a_{ijst}f(\x_1,\cdots,\widehat{\x_i},\cdots, \widehat{\x_s},\cdots,\xst, \cdots, \xij,\cdots,\x_{m+2},z)\tag{a4}\\
     &-\sum_{i<s<j<t}\widetilde{a_{ijst}}f(\x_1,\cdots,\widehat{\x_i},\cdots, \widehat{\x_s},\cdots,\xij, \cdots, \xst,\cdots,\x_{m+2},z)\tag{a5}\\
     &-\sum_{i<j<s<t}a_{ijst}f(\x_1,\cdots,\widehat{\x_i},\cdots,\xij,\cdots, \widehat{\x_s}, \cdots, \xst,\cdots, \x_{m+2},z)\tag{a6}\\
 &+\sum_{k<i<j}\widetilde{b_{ijk}}f(\x_1,\cdots, \widehat{\x_k},\cdots, \widehat{\x_i}, \cdots, \x_k\cdot(\xij), \cdots, \x_{m+2},z)\tag{b1}\\
  &-\sum_{i<k<j}b_{ijk}f(\x_1,\cdots, \widehat{\x_i},\cdots, \widehat{\x_k}, \cdots, \x_k\cdot(\xij), \cdots, \x_{m+2},z)\tag{b2}\\
   &-\sum_{i<j<k}\widetilde{b_{ikj}}f(\x_1,\cdots, \widehat{\x_i},\cdots, \widehat{\x_j}, \cdots, (\xij)\cdot\x_k, \cdots, \x_{m+2},z)\tag{b3}\\
 &+\sum_{k<i<j}c_{ijk}f(\x_1,\cdots,\widehat{\x_k},\cdots,\widehat{\x_i},\cdots,\xij,\cdots,\x_{m+2}, \x_k\cdot z)\tag{c1}\\
  &-\sum_{i<k<j}\widetilde{c_{ijk}}f(\x_1,\cdots,\widehat{\x_i},\cdots,\widehat{\x_k},\cdots,\xij,\cdots,\x_{m+2}, \x_k\cdot z)\tag{c2}\\
   &-\sum_{i<j<k}c_{ijk}f(\x_1,\cdots,\widehat{\x_i},\cdots,\xij,\cdots,\widehat{\x_k},\cdots,\x_{m+2}, \x_k\cdot z)\tag{c3}\\
 &-\sum_{i<j}\widetilde{d_{ij}}f(\x_1,\cdots,\widehat{\x_i},\cdots,\widehat{\x_j},\cdots,\x_{m+2}, (\xij)\cdot z)\tag{d1}\\
 &+\sum_{k<i<j}e_{ijk}\x_k\cdot f(\x_1,\cdots,\widehat{\x_k},\cdots,\widehat{\x_i},\cdots,\xij,\cdots,\x_{m+2}, z)\tag{e1}\\
  &-\sum_{i<k<j}\widetilde{e_{ijk}}\x_k\cdot f(\x_1,\cdots,\widehat{\x_i},\cdots,\widehat{\x_k},\cdots,\xij,\cdots,\x_{m+2}, z)\tag{e2}\\
   &-\sum_{i<j<k}e_{ijk}\x_k\cdot f(\x_1,\cdots,\widehat{\x_i},\cdots,\xij,\cdots,\widehat{\x_k},\cdots,\x_{m+2}, z)\tag{e3}\\
 &-\sum_{i<j}\widetilde{g_{ij}}(\xij)\cdot f(\x_1,\cdots,\widehat{\x_i},\cdots,\widehat{\x_j},\cdots,\x_{m+2},z)\tag{g1}\\
 &+\sum_{i<j\leq {m+1}}h_{ij}(f(\x_1,\cdots,\widehat{\x_i},\cdots,\xij,\cdots,\x_{m+1},~~)\cdot\x_{m+2})\cdot z\tag{h1}\\
 \begin{split}
 &+\sum_{k=1}^{m+1}(-1)^{k+m}(-1)^{|\x_k|(|\x_{k+1}|+\cdots+|\x_{m+1}|)}\\
 &\relphantom{+\sum_{k=1}^{m+1}}\cdot(f(\x_1,\cdots,\widehat{\x_k},\cdots,\x_{m+1},~~)\cdot(\x_k\cdot\x_{m+2}))\cdot z
 \end{split}\tag{l1}\\
 &+\sum_{s<t<i}c_{sti}f(\x_1,\cdots,\widehat{\x_s},\cdots,\xst,\cdots,\widehat{\x_i},\cdots,\x_{m+2}, \x_i\cdot z)\tag{c4}\\
  &+\sum_{s<i<t}\widetilde{c_{sti}}f(\x_1,\cdots,\widehat{\x_s},\cdots,\widehat{\x_i},\cdots,\xst,\cdots,\x_{m+2}, \x_i\cdot z)\tag{c5}\\
   &-\sum_{i<s<t}c_{sti}f(\x_1,\cdots,\widehat{\x_i},\cdots,\widehat{\x_s},\cdots,\xst,\cdots,\x_{m+2}, \x_i\cdot z)\tag{c6}\\
 &+\sum_{k<i}\widetilde{d_{ik}}f(\x_1,\cdots,\widehat{\x_k},\cdots,\widehat{\x_i},\cdots,\x_{m+2}, \x_k\cdot(\x_i\cdot z))\tag{d2}\\
  &-\sum_{i<k}d_{ik}f(\x_1,\cdots,\widehat{\x_i},\cdots,\widehat{\x_k},\cdots,\x_{m+2}, \x_k\cdot(\x_i\cdot z))\tag{d3}\\
 &+\sum_{k<i}p_{ki}\x_k\cdot f(\x_1,\cdots,\widehat{\x_k},\cdots,\widehat{\x_i},\cdots,\x_{m+2}, \x_i\cdot z)\tag{p1}\\
  &-\sum_{i<k}\widetilde{p_{ki}}\x_k\cdot f(\x_1,\cdots,\widehat{\x_i},\cdots,\widehat{\x_k},\cdots,\x_{m+2}, \x_i\cdot z)\tag{p2}\\
  \begin{split}
 &+\sum_{i=1}^{m+1} (-1)^{i+m}(-1)^{|\x_i|(|\x_{i+1}|+\cdots+|\x_{m+2}|)}\\
 &\relphantom{-\sum_{i=1}^{n-1}}\cdot(f(\x_1,\cdots,\widehat{\x_i},\cdots,\x_{m+1},~~)\cdot\x_{m+2})\cdot(\x_i\cdot z)
 \end{split}\tag{l2}\\
 &+(f(\x_1,\cdots,\x_{m},~~)\cdot \x_{m+1})\cdot(\x_{m+2}\cdot z)\tag{q1}\\
 &+\sum_{s<t<i}e_{sti}\x_i\cdot f(\x_1,\cdots,\widehat{\x_s},\cdots,\xst,\cdots,\widehat{\x_i},\cdots,\x_{m+2}, z)\tag{e4}\\
  &+\sum_{s<i<t}\widetilde{e_{sti}}\x_i\cdot f(\x_1,\cdots,\widehat{\x_s},\cdots,\widehat{\x_i},\cdots,\xst,\cdots,\x_{m+2}, z)\tag{e5}\\
   &-\sum_{i<s<t}e_{sti}\x_i\cdot f(\x_1,\cdots,\widehat{\x_i},\cdots,\widehat{\x_s},\cdots,\xst,\cdots,\x_{m+2}, z)\tag{e6}\\
   &+\sum_{k<i}\widetilde{p_{ik}}\x_i\cdot f(\x_1,\cdots,\widehat{\x_k},\cdots,\widehat{\x_i},\cdots,\x_{m+2}, \x_k\cdot z)\tag{p3}\\
  &-\sum_{i<k}p_{ik}\x_i\cdot f(\x_1,\cdots,\widehat{\x_i},\cdots,\widehat{\x_k},\cdots,\x_{m+2}, \x_k\cdot z)\tag{p4}\\
 &-\sum_{k<i}g_{ki}\x_k\cdot(\x_i\cdot f(\x_1,\cdots,\widehat{\x_k},\cdots,\widehat{\x_i},\cdots,\x_{m+2},z))\tag{g2}\\
  &+\sum_{i<k}\widetilde{g_{ki}}\x_k\cdot(\x_i\cdot  f(\x_1,\cdots,\widehat{\x_i},\cdots,\widehat{\x_k},\cdots,\x_{m+2},z))\tag{g3}\\
  \begin{split}
 &-\sum_{i=1}^{m+1} (-1)^{i+m}(-1)^{|\x_i|(|f|+|\x_{1}|+\cdots+|\x_{i-1}|)}\\
 &\relphantom{-\sum_{i=1}^{n-1}}\x_i\cdot((f(\x_1,\cdots,\widehat{\x_i},\cdots,\x_{m+1},~~)\cdot \x_{m+2})\cdot z)
 \end{split}\tag{l3}\\
 &-(-1)^{|\x_{m+2}|(|f|+|\x_1|+\cdots+|\x_{m+1}|)}\x_{m+2}\cdot((f(\x_1,\cdots,\x_{m},~~)\cdot \x_{m+1})\cdot z)\tag{q2}\\
 &-\sum_{s<t\leq {m+1}}h_{st}(f(\x_1,\cdots,\widehat{\x_s},\cdots,\xst,\cdots,\x_{m+1},~~)\cdot\x_{m+2})\cdot z\tag{h2}\\
 \begin{split}
 &-\!\!\sum_{i=1}^{n-1}\sum_{k=1}^{m+1} (-1)^{m+k}(-1)^{(|f|+|\x_{1}|+\cdots+|\x_{{m+1}}|)(|\x_{m+2}^1|+\cdots+|\x_{m+2}^{i-1}|)} (-1)^{|\x_k|(|\x_{k+1}|+\cdots+|\x_{m+1}|)}\\
 &\relphantom{-\sum_{i=1}^{n-1}\sum_{k=1}^{m+1}}\cdot[\x_{m+2}^1,\cdots,f(\x_1,\cdots,\widehat{\x_k},\cdots,\x_{m+1}, \x_k\cdot\x_{m+2}^i),\cdots,\x_{m+2}^{n-1},z]
 \end{split}\tag{l4}\\
 \begin{split}
 &+\!\!\sum_{i=1}^{n-1}\sum_{k=1}^{m+1} (-1)^{m+k}(-1)^{(|f|+|\x_{1}|+\!\cdots+|\x_{{m+1}}|)(|\x_{m+2}^1|+\cdots+|\x_{m+2}^{i-1}|)} (-1)^{|\x_k|(|f|+|\x_1|+\cdots+|\x_{k-1}|)}\\
 &\relphantom{+\!\sum_{i=1}^{n-1}\sum_{k=1}^{m+1}}\cdot[\x_{m+2}^1,\cdots,\x_k\cdot f(\x_1,\cdots,\widehat{\x_k},\cdots,\x_{m+1},\x_{m+2}^i),\cdots,\x_{m+2}^{n-1},z]
 \end{split}\tag{l5}\\
 \begin{split}
 &-\sum_{i=1}^{n-1}(-1)^{(|f|+|\x_{1}|+\cdots+|\x_{{m+1}}|)(|\x_{m+2}^1|+\cdots+|\x_{m+2}^{i-1}|)}\notag\\
 &\relphantom{-\sum_{i=1}^{n-1}}\cdot[\x_{m+2}^1,\cdots,(f(\x_1,\cdots,\x_{m},~~)\cdot\x_{m+1})\cdot\x_{m+2}^i,\cdots,\x_{m+2}^{n-1},z],
 \end{split}\tag{q3}
\end{align}
where
\begin{align*}
a_{ijst}=&(-1)^{s+i}(-1)^{|\x_i|(|\x_{i+1}|+\cdots+|\x_{j-1}|)}(-1)^{|\x_s|(|\x_{s+1}|+\cdots+|\x_{t-1}|)},
& \widetilde{a_{ijst}}=&(-1)^{|\x_i||\x_s|}a_{ijst};\\
b_{ijk}=&(-1)^{i+k}(-1)^{|\x_i|(|\x_{i+1}|+\cdots+|\x_{j-1}|)}(-1)^{|\x_k|(|\x_{k+1}|+\cdots+|\x_{j-1}|)},
& \widetilde{b_{ijk}}=&(-1)^{|\x_i||\x_k|}b_{ijk};\\
c_{ijk}=&(-1)^{i+k}(-1)^{|\x_i|(|\x_{i+1}|+\cdots+|\x_{j-1}|)}(-1)^{|\x_k|(|\x_{k+1}|+\cdots+|\x_{m+2}|)},
& \widetilde{c_{ijk}}=&(-1)^{|\x_i||\x_k|}c_{ijk};\\
d_{ij}=&(-1)^{i+j}(-1)^{|\x_i|(|\x_{i+1}|+\cdots+|\x_{m+2}|)}(-1)^{|\x_j|(|\x_{j+1}|+\cdots+|\x_{m+2}|)},
& \widetilde{d_{ij}}=&(-1)^{|\x_i||\x_j|}d_{ij};\\
e_{ijk}=&(-1)^{i+k+1}(-1)^{|\x_i|(|\x_{i+1}|+\cdots+|\x_{j-1}|)}(-1)^{|\x_k|(|f|+|\x_{1}|+\cdots+|\x_{k-1}|)},
& \widetilde{e_{ijk}}=&(-1)^{|\x_i||\x_k|}e_{ijk};\\
g_{ij}=&(-1)^{i+j+1}(-1)^{|\x_i|(|f|+|\x_{1}|+\cdots+|\x_{i-1}|)}(-1)^{|\x_j|(|f|+|\x_{1}|+\cdots+|\x_{j-1}|)},
& \widetilde{g_{ij}}=&(-1)^{|\x_i||\x_j|}g_{ij};\\
h_{ij}=&(-1)^{i+m}(-1)^{|\x_i|(|\x_{i+1}|+\cdots+|\x_{j-1}|)},
& \widetilde{h_{ij}}=&(-1)^{|\x_i||\x_j|}h_{ij};\\
%l_k=&(-1)^{k+m}(-1)^{|\x_k|(|\x_{k+1}|+\cdots+|\x_{m+1}|)},
%& \widetilde{l_k}=&(-1)^{|\x_k||\x_{m+2}|}l_k;\\
p_{ki}=&(-1)^{i+k+1}(-1)^{|\x_i|(|\x_{i+1}|+\cdots+|\x_{m+2}|)}(-1)^{|\x_k|(|f|+|\x_{1}|+\cdots+|\x_{k-1}|)},
&\widetilde{p_{ki}}=&(-1)^{|\x_i||\x_k|}p_{ki}.
\end{align*}

It can be verified that the sum of terms labeled with the same letter vanishes(e.g. (a1)+$\cdots$+(a6)=0), then $\delta^2=0$ and $\delta$ is called a coboundary operator. Therefore, we get the following theorem.

\bthm
The coboundary operator $\delta$ introduced in Definition \ref{def:delta} satisfies $\delta^2f=0, \forall f\in C^{m}(\g, V)$.
\ethm

\bre
The coboundary operator $\delta$ as above is a generalization of which of $n$-Lie algebras in \cite{AI} and of Lie superalgebras in \cite{SZ}.
\ere

The map $f\in C^{m}(\g, V)$ is called an $m$-supercocycle if $\delta f=0$. We denote by $Z^{m}(\g,V)$ the graded subspace spanned by $m$-supercocycles. Since $\delta^{2}f=0$ for all $f \in C^{m}(\g, V)$, $\delta C^{m-1}(\g, V)$ is a graded subspace of $Z^{m}(\g,V)$. Therefore we can define a graded cohomology space $H^{m}(\g,V)$ of $\g$ as the graded factor space $Z^{m}(\g,V)/\delta C^{m-1}(\g, V).$
\section{Extension of first-class $n$-Lie superalgebras}

Let $\g, \mathfrak{a}, \mathfrak{b}$ be first-class $n$-Lie superalgebras over $\K$. $\g$ is called an extension of $\mathfrak{b}$ by $\mathfrak{a}$ if there is an exact sequence of first-class $n$-Lie superalgebras:
$$\xymatrix{0\ar[r]& \mathfrak{a}\ar[r]^\iota& \g\ar[r]^\pi& \mathfrak{b}\ar[r]& 0}.$$

Suppose that $\mathfrak{a}$ is an abelian graded ideal of $\g$, i.e., $\mathfrak{a}$ is a graded ideal such that $[\mathfrak{a},\mathfrak{a},\underbrace{\g,\cdots,\g}_{n-2}]=0$. We consider the case that $\g$ is an extension of $\mathfrak{b}$ by an abelian graded ideal $\mathfrak{a}$ of $\g$. Let $\tau:\mathfrak{b}\rightarrow \g$ be a homogeneous linear map of degree 0 with $\pi\circ\tau=\id_{\mathfrak{b}}$.
Let $\B=b_1\wedge\cdots\wedge b_{n-1}\in\mathfrak{b}^{\wedge^{n-1}}$ and let $\rho: \mathfrak{b}^{\wedge^{n-1}}\rightarrow \End(\mathfrak{a}), \B\mapsto\tau(\B)=\tau(b_1)\wedge\cdots\wedge\tau(b_{n-1})$. Then $\mathfrak{a}$ becomes a graded $\mathfrak{b}$-module. Let us write $\tau(b)=(0,b)$ and then denote the elements of $\g$ by $(a, b)$ for all $a\in \mathfrak{a}$ and $b\in\mathfrak{b}$. Then, the bracket in $\g$ is defined by
\beq\label{eq:n-def}
[(a_1,b_1),\cdots,(a_n,b_n)]
=\left(\sum_{i=1}^n[\tau(b_1),\cdots,a_i,\cdots,\tau(b_n)]+f(\B,b_n), ~\B\cdot b_n\right),
\eeq
where  $f(\B,b_n)=\tau(\B)\cdot\tau(b_n)-\tau(\B\cdot b_n)$ and  $|(a_i,b_i)|=|a_i|=|b_i|, \forall 1\leq i\leq n$. Then $f\in C^1(\mathfrak{b},\mathfrak{a})_{\bar{0}}$.
Let $\A=a_1\wedge\cdots\wedge a_{n-1}$ and $(\A,\B)=(a_1,b_1)\wedge\cdots\wedge(a_{n-1},b_{n-1})$. Then
\begin{align*}
&(\A,\B)\cdot((\A',\B')\cdot(a_n',b_n'))\\
&-\sum_{i=1}^n(-1)^{|\A|(|a_1'|+\cdots+|a_{i-1}'|)}[(a_1',b_1'),\cdots,(\A,\B)\cdot(a_i',b_i'),\cdots,(a_n',b_n')]\\
=&(\A,\B)\cdot\left(\sum_{i=1}^n[\tau(b_1'),\cdots,a_i',\cdots,\tau(b_n')]+f(\B', b_n'), ~\B'\cdot b_n'\right)\\
&-\sum_{i=1}^n(-1)^{|\A|(|a_1'|+\cdots+|a_{i-1}'|)}\\
&\cdot\left[(a_1',b_1'),\cdots,
\left(\left\{\begin{aligned}
&\sum_{j=1}^{n-1}[\tau(b_1),\cdots, a_j,\cdots,\tau(b_{n-1}),\tau(b_i')]\\
&+\tau(\B)\cdot a_i'+f(\B,b_i')
\end{aligned}\right\}
, ~\B\cdot b_i'\right)
,\cdots,(a_n',b_n')\right]\\
=&\left(\left\{\begin{aligned}
&\tau(\B)\cdot\left(\sum_{i=1}^n[\tau(b_1'),\cdots,a_i',\cdots,\tau(b_n')]\right)+\tau(\B)\cdot f(\B',b_n')\\
&+\sum_{j=1}^{n-1}\left[\tau(b_1),\cdots, a_j,\cdots,\tau(b_{n-1}),\tau(\B'\cdot b_n')\right]+f(\B, \B'\cdot b_n')
\end{aligned}\right\}
, ~\B\cdot(\B'\cdot b_n)\right)\\
&-\sum_{i=1}^n(-1)^{|\A|(|a_1'|+\cdots+|a_{i-1}'|)}\\
&\cdot\left(\left\{\begin{aligned}
&\sum_{j=1}^{n-1}\Big[\tau(b_1'),\cdots,[\tau(b_1),\cdots, a_j,\cdots,\tau(b_{n-1}),\tau(b_i')],\\
&\relphantom{\sum_{j=1}^{n-1}\Big[}\cdots,\tau(b_n')\Big]\\
&+[\tau(b_1'),\cdots,\tau(\B)\cdot a_i',\cdots,\tau(b_n')]\\
&+[\tau(b_1'),\cdots,f(\B,b_i'),\cdots,\tau(b_n')]\\
&+\sum_{j\neq i}[\tau(b_1'),\cdots,a_j',\cdots,\tau(\B\cdot b_i'),\cdots,\tau(b_n')]\\
&+f(b_1',\cdots,\B\cdot b_i',\cdots,b_n')
\end{aligned}\right\}
, ~[b_1',\cdots,\B\cdot b_i',\cdots,b_n']\right)\\
%=&\left(\begin{aligned}
%&\sum_{i=1}^nt_i^{(1)}\tau(\B)\cdot\left([\tau(b_1'),\cdots,\widehat{\tau(b_i')},\cdots,\tau(b_n'),a_i']+f(\B',b_n')\right)\\
%&+\sum_{j=1}^nt_j^{(2)}\left[\tau(b_1),\cdots,\widehat{\tau(b_j)},\cdots,\tau(b_{n-1}),\tau(\B'\cdot b_n'), a_j\right]
%\end{aligned}
%,~\B\cdot(\B'\cdot b_n)\right)\\
%&-\left(\begin{aligned}
%&\sum_{i=1}^n(-1)^{n-i}t_i^{(3)}[\tau(b_1'),\cdots,\widehat{\tau(b_i')},\cdots,\tau(b_n'),\tau(\B)\cdot a_i']\\
%&t_i^{(4)}[\tau(b_1'),\cdots,\widehat{\tau(b_i')},\cdots,\tau(b_n'),[\tau(b_1),\cdots,\widehat{\tau(b_j)},\cdots,\tau(b_{n-1}),\tau(b_i'), a_j]]\\
%&+t_i^{(5)}[\tau(b_1'),\cdots,\widehat{\tau(b_i')},\cdots,\tau(b_n'), f(b_1,\cdots,b_{n-1},b_i')]\\
%&+t_i^{(6)} [\tau(b_1'),\cdots,\widehat{\tau(b_j')},\cdots,\tau(\B'\cdot b_n'),\cdots,\tau(b_n'),a_j']\\
%&+t_i^{(7)} [\tau(b_1'),\cdots,\tau(\B'\cdot b_n'),\cdots,\widehat{\tau(b_j')},\cdots,\tau(b_n'),a_j']\\
%&+f(\B\cdot\B', b_n')+(-1)^{|\A||\A'|}f(\B,\B'\cdot b_n')
%\end{aligned}
%~,\B\cdot(\B'\cdot b_n)\right)\\
=&(\delta f(\B,\B',b_n'), 0).
\end{align*}
Therefore, $f\in Z^1(\mathfrak{b}, \mathfrak{a})_{\bar{0}}$.

Conversely, suppose that an abelian first-class $n$-Lie superalgebra $\mathfrak{a}$ is a graded $\mathfrak{b}$-module, $\rho(\B)\cdot a:=\tau(\B)\cdot a$, and $f\in Z^1(\mathfrak{b}, \mathfrak{a})_{\bar{0}}$.
Let $\g=\mathfrak{a}\times\mathfrak{b}$. Then $\g$ is a first-class $n$-Lie superalgebra with the bracket defined by (\ref{eq:n-def}). Then we can define an exact sequence
$$\xymatrix{0\ar[r]& \mathfrak{a}\ar[r]^\iota& \g\ar[r]^\pi& \mathfrak{b}\ar[r]& 0},$$
where $\iota(a)=(a,0), \pi(a,b)=b$. Thus $\g$ is an extension of $\mathfrak{b}$ by $\mathfrak{a}$ and $\iota(\mathfrak{a})$ is an abelian graded ideal of $\g$.

Therefore, we get the following theorem.
\bthm
Suppose that $\mathfrak{a}, \mathfrak{b}$ are first-class $n$-Lie superalgebras over $\K$ and $\mathfrak{a}$ is abelian. Then there is a one-to-one correspondence between extensions of $\mathfrak{b}$ by $\mathfrak{a}$ and $Z^1(\mathfrak{b}, \mathfrak{a})_{\bar{0}}$.
\ethm
\section{$T$*-extension of first-class $n$-Lie superalgebras}

Let $\g$ be a first-class $n$-Lie superalgebra, $\g^*$ be its dual space.
Since $\g=\g_{\bar{0}}\oplus\g_{\bar{1}}$ and $\g^*=\g^*_{\bar{0}}\oplus\g^*_{\bar{1}}$ are $\Z_2$-graded vector space, the direct sum $\g\oplus\g^*=(\g_{\bar{0}}\oplus\g^*_{\bar{0}})
\oplus(\g_{\bar{1}}\oplus\g^*_{\bar{1}})$ is a $\Z_2$-graded vector space. In the sequel, whenever $x+f\in \g\oplus\g^*$ appears, it means that $x+f$ is homogeneous and $|x+f|=|x|=|f|$.

Let $\theta$ be a homogeneous $n$-linear map from $\g^{\wedge^n}$ into $\g^*$ of degree 0. Now we define a bracket on $\g\oplus\g^*$:
\beq\label{eq:bracketofTextension}\begin{split}
[x_1+f_1,\cdots,x_n+f_n]_{\theta}=&[x_1,\cdots,x_n]_{\g}+\theta(x_1,\cdots,x_n)\\
&+ \sum_{i=1}^n(-1)^{n-i}(-1)^{|x_i|(|x_{i+1}|+\cdots+|x_n|)} \ad^*(x_1,\cdots,\widehat{x_i},\cdots,x_n)\cdot f_i.
\end{split}\eeq
\blem
$\g\oplus\g^*$ is a first-class $n$-Lie superalgebra if and only if $\theta\in Z^1(\g, \g^*)_{\bar{0}}$.
\elem
\bpf
It's clear that $[ ,\cdots, ]_{\theta}$ satisfies (\ref{eq:2}) if and only if $\theta\in C^1(\g, \g^*)_{\bar{0}}$. Let $\x+\F=(x_1+f_1)\wedge\cdots\wedge(x_{n-1}+f_{n-1})$ and $\y+\G=(y_1+g_1)\wedge\cdots\wedge(y_{n-1}+g_{n-1})$. Then we have
\begin{align*}
&(\x+\F)\cdot \left((\y+\G)\cdot (y_n+g_n)\right)\\
=&(\x+\F)\cdot\Big\{
\sum_{i=1}^n(-1)^{n-i}(-1)^{|y_i|(|y_{i+1}|+\cdots+|y_n|)} \ad^*(y_1,\cdots,\widehat{y_i},\cdots,y_n)\cdot g_i\\
&\relphantom{(\x+\F)\cdot\Big\{}+\y\cdot y_n+\theta(\y,y_n)\Big\}\\
=&\x\cdot(\y\cdot y_n)+\theta(\x,\y\cdot y_n)+\ad^*(\x)\cdot\theta(\y,y_n)\\
&+\sum_{j=1}^{n-1}(-1)^{n-j}(-1)^{|x_j|(|x_{j+1}|+\cdots+|x_{n-1}|+|\y|+|y_n|)} \ad^*(x_1,\cdots,\widehat{x_j},\cdots,x_{n-1},\y\cdot y_n)\cdot f_j\\
&+\sum_{i=1}^n(-1)^{n-i}(-1)^{|y_i|(|y_{i+1}|+\cdots+|y_n|)} \ad^*(\x)\cdot(\ad^*(y_1,\cdots,\widehat{y_i},\cdots,y_n)\cdot g_i)
\end{align*}
and
\begin{align*}
&\sum_{i=1}^n(-1)^{|\x|(|y_1|+\cdots+|y_{i-1}|)}[y_1+g_1,\cdots,(\x+\F)\cdot (y_i+g_i),\cdots,y_n+g_n]_{\theta}\\
=&\sum_{i=1}^n(-1)^{|\x|(|y_1|+\cdots+|y_{i-1}|)}\bigg[y_1+g_1,\cdots, \Big\{\x\cdot y_i+\theta(\x,y_i)+\ad^*(\x)\cdot g_i\\
&+\sum_{j=1}^{n-1}(-1)^{n-j}(-1)^{|x_j|(|x_{j+1}|+\cdots+|x_{n-1}|+|y_i|)}
\ad^*(x_1,\!\cdots\!,\widehat{x_j},\!\cdots\!,x_{n-1},y_i)\cdot f_j\Big\},\!\cdots\!,y_n+g_n\bigg]_\theta\\
=&\sum_{i=1}^n(-1)^{|\x|(|y_1|+\cdots+|y_{i-1}|)}\bigg\{[y_1,\cdots,\x\cdot y_i,\cdots,y_n]+\theta(y_1,\cdots,\x\cdot y_i,\cdots,y_n)\\
&\relphantom{+}+\sum_{k<i}(-1)^{n-k}(-1)^{|y_k|(|y_{k+1}|+\cdots+|y_n|+|\x|)}
\ad^*(y_1,\cdots,\widehat{y_k},\cdots,\x\cdot y_i,\cdots,y_n)\cdot g_k\\
&\relphantom{+}+\sum_{i<k}(-1)^{n-k}(-1)^{|y_k|(|y_{k+1}|+\cdots+|y_n|)}
\ad^*(y_1,\cdots,\x\cdot y_i,\cdots,\widehat{y_k},\cdots,y_n)\cdot g_k\\
&\relphantom{+}+(-1)^{n-i}(-1)^{(|\x|+|y_i|)(|y_{i+1}|+\cdots+|y_n|)}
\ad^*(y_1,\cdots,\widehat{y_i},\cdots,y_n)\cdot \Big\{
\theta(\x,y_i)+\ad^*(\x)\cdot g_i\\
&\relphantom{+}\relphantom{+}+\sum_{j=1}^{n-1}(-1)^{n-j}(-1)^{|x_j|(|x_{j+1}|+\cdots+|x_{n-1}|+|y_i|)}
\ad^*(x_1,\cdots,\widehat{x_j},\cdots,x_{n-1},y_i)\cdot f_j\Big\}\bigg\}.
\end{align*}
Since $[ ,\cdots, ]_{\g}$ satisfies (\ref{eq:3'}) and $\ad^*(\x)$ satisfies (\ref{eq:rho3}), it can be concluded that $[ ,\cdots, ]_{\theta}$ satisfies (\ref{eq:3'}) if and only if
\begin{align*}
0=&\theta(\x,\y\cdot y_n)+\ad^*(\x)\cdot\theta(\y,y_n)
  -\sum_{i=1}^n(-1)^{|\x|(|y_1|+\cdots+|y_{i-1}|)}\theta(y_1,\cdots,\x\cdot y_i,\cdots,y_n)\\
  &-\sum_{i=1}^n(-1)^{|\x|(|y_1|+\cdots+|y_{i-1}|)}
  (-1)^{n-i}(-1)^{(|\x|+|y_i|)(|y_{i+1}|+\cdots+|y_n|)}\\
  &\relphantom{+}\relphantom{+} \cdot \ad^*(y_1,\cdots,\widehat{y_i},\cdots,y_n)\cdot\theta(\x,y_i)\\
 =&\delta\theta(\x,\y,y_n),
\end{align*}
i.e., $\theta\in Z^1(\g, \g^*)_{\bar{0}}$.
\epf

\bdefn
Let $\g$ be a first-class $n$-Lie superalgebra. A bilinear form $\langle ,\rangle_{\g}$ on $\g$ is said to be nondegenerate if
$$\g^\perp=\{x\in \g|\langle x,y\rangle_{\g}=0, \forall y\in \g\}=0;$$
invariant if
$$\langle[x_1, \!\cdots\!, x_{n-1},y]_{\g},z\rangle_{\g}=-(-1)^{(|x_1|+\cdots+|x_{n-1}|)|y|}\langle y,[x_1, \!\cdots\!, x_{n-1},z]_{\g}\rangle_{\g}, \forall x_1, \!\cdots\!, x_{n-1}, y, z\in \g;$$
supersymmetric if
$$\langle x,y\rangle_{\g}=(-1)^{|x||y|}\langle y,x\rangle_{\g};$$
consistent if
$$\langle x,y\rangle_{\g}=0, \forall x, y\in \g, |x|\neq|y|.$$
In this section, we only consider consistent bilinear forms. If $\g$ admits a nondegenerate invariant supersymmetric bilinear form $\langle ,\rangle_{\g}$, then we call $(\g, \langle ,\rangle_{\g})$ a metric first-class $n$-Lie superalgebra.
\edefn

\blem
Define a bilinear form $\langle , \rangle_{\theta}:(\g\oplus\g^*)\times (\g\oplus\g^*)\rightarrow \K$ by
$$\langle x+f, y+g \rangle_{\theta}=f(y)+(-1)^{|x||y|}g(x).$$ Then $\langle y+g, x+f \rangle_{\theta}=(-1)^{|x||y|}\langle x+f, y+g \rangle_{\theta}$ and  $\langle , \rangle_{\theta}$ is nondegenerate. Moreover, $(\g\oplus\g^*, \langle , \rangle_{\theta})$ is metric if and only if the following identity holds:
\beq\label{eq:supercyclic} \theta(\x,y)(z)+(-1)^{|y||z|}\theta(\x,z)(y)=0.\eeq
\elem
\bpf
$(\g\oplus\g^*, \langle , \rangle_{\theta})$ is metric if and only if
\begin{align*}
0=&\langle(\x+\F)\cdot(y+g),z+h\rangle_{\theta}
+(-1)^{|\x||y|}\langle y+g,(\x+\F)\cdot(z+h)\rangle_{\theta}\\
=&\langle\x\cdot y+\theta(\x,y)+\ad^*(\x)\cdot g, z+h\rangle_{\theta}\\ &+\left\langle\sum_{i=1}^{n-1}(-1)^{n-i}(-1)^{|x_i|(|x_{i+1}|+\cdots+|x_{n-1}|+|y|)}\ad^*\!(x_1,\cdots,\widehat{x_i},\cdots,x_{n-1},y)\cdot f_i, z+h\right\rangle_{\theta}\\
&+(-1)^{|\x||y|}\left\langle y+g, \x\cdot z+\theta(\x, z)+ \ad^*(\x)\cdot h\right\rangle_{\theta}\\
&+(-1)^{|\x||y|}\!\left\langle \!y+g, \!\sum_{i=1}^{n-1}(-1)^{n-i}(-1)^{|x_i|(|x_{i+1}|+\cdots+|x_{n-1}|+|z|)}\ad^*(x_1,\!\cdots\!,\widehat{x_i},\!\cdots\!,x_{n-1},z)\!\cdot\! f_i\!\right\rangle_{\theta}\\
=&\theta(\x,y)(z)+(-1)^{|y||z|}\theta(\x,z)(y),
\end{align*}
i.e., (\ref{eq:supercyclic}) holds.\epf

Now we give the definition of $T^*$-extensions.
\bdefn
For a 1-supercocycle $\theta$ satisfying (\ref{eq:supercyclic}) we shall call the metric first-class $n$-Lie superalgebra $(\g\oplus \g^{*},\langle , \rangle_{\theta})$ the $T^*$-extension of $\g$ (by $\theta$) and denote it by $T_\theta^*\g$.
\edefn

\bthm
Let $\g$ be a first-class $n$-Lie superalgebra over a field $\K$. Let
\Beq \g^{(0)}=\g, \g^{(m+1)}=[\g^{(m)},\cdots,\g^{(m)}]_{\g} \text{~~and~~} \g^{0}=\g, \g^{m+1}=[\g^{m},\g,\cdots,\g]_{\g}, \forall m\geq0. \Eeq
$\g$ is called solvable (nilpotent) of length $k$ if and only if there is a smallest integer $k$ such that $\g^{(k)}=0$ ($\g^{k}=0$). Then
\begin{enumerate}[(1)]
   \item  If $\g$ is solvable of length $k$, then $T^{*}_{\theta}\g$ is solvable of length $k$ or $k+1$.
   \item  If $\g$ is nilpotent of length $k$, then $T^{*}_{\theta}\g$ is nilpotent of length at least $k$ and at most $2k-1$. In particular, the nilpotent length of $T^{*}_{0}\g$ is $k$.
   \item  If $\g$ can be decomposed into a direct sum of two (graded) ideals of $\g$, then $T^{*}_{0}\g$ can be too.
\end{enumerate}
\ethm
\bpf
(1) Suppose that $\g$ is solvable of length $k$. Since
$(T^{*}_{\theta}\g)^{(m)}/\g^{*}\cong \g^{(m)}$ and $\g^{(k)}=0$, we have
 $(T^{*}_{\theta}\g)^{(k)}\subseteq \g^{*}$, which implies $(T^{*}_{\theta}\g)^{(k+1)}=0$ because $\g^{*}$ is abelian, and it follows that $T^{*}_{\theta}\g$ is
solvable of length $k$ or $k+1$.

(2) Suppose that $\g$ is nilpotent of length $k$. Since $(T^{*}_{\theta}\g)^{m}/\g^{*}\cong \g^{m}$ and $\g^{k}=0$, we have
$(T^{*}_{\theta}\g)^{k}\subseteq \g^{*}$. Let $f\in(T^{*}_{\theta}\g)^{k}\subseteq \g^{*}, y\in \g$, $\x_{j}+\F_j=(\x_{j}^1+\F_j^1)\wedge\cdots\wedge(\x_{j}^{n-1}+\F_j^{n-1})\in (T^{*}_{\theta}\g)^{\wedge^{n-1}}$, $j=1,\cdots,k-1$. Then
$$\left((\x_{1}+\F_1)\cdots(\x_{k-1}+\F_{k-1})\cdot f\right)(y)
=(\ad^*(\x_1)\cdots \ad^*(\x_{k-1})\cdot f)(y)\in f(\g^k)=0.$$
This proves that $(T^{*}_{\theta}\g)^{2k-1}=0$. Hence $T^{*}_{\theta}\g$ is nilpotent of length at least $k$ and at most $2k-1$.

Now consider the case of trivial $T^*$-extension $T^{*}_{0}\g$ of $\g$. Note that
\begin{align*}
&(\x_{1}+\F_1)\cdots(\x_{k-1}+\F_{k-1})\cdot (y+g)\\
=&\ad(\x_1)\cdots \ad(\x_{k-1})\cdot y+\ad^*(\x_1)\cdots \ad^*(\x_{k-1})\cdot g\\
&+\sum_{j=1}^{k-1}\sum_{i=1}^{n-1}
(-1)^{n-i}(-1)^{|\x_j^i|(|\x_j^{i+1}|+\cdots+|\x_j^{n-1}|+|y|+|\x_{j+1}|+\cdots+|\x_{k-1}|)}\\
&\relphantom{+}\cdot\ad^*(\x_1)\cdots \ad^*(\x_{j-1})\ad^*(\x_j^1,\cdots,\widehat{\x_j^i},\cdots,\x_j^{n-1}, \ad(\x_{j+1})\cdots\ad(\x_{k-1})\cdot y)\cdot\F_j^i\\
=&0.
\end{align*}
Then $(T^{*}_{\theta}\g)^k=0$, as required.

(3) Suppose that $0\neq \g=I\oplus J$,  where $I$ and $J$ are two nonzero (graded) ideals of $\g$. Let $I^{*}=\{f\in\g^*| f(J)=0\}$ and $J^{*}=\{f\in\g^*| f(I)=0\}$. Then $I^{*}$(resp. $J^{*}$) can canonically be identified with the dual space of $I$(resp. $J$) and $\g^*\cong I^*\oplus J^*$.

Note that
\begin{align*}
[T^{*}_{0}I,T^{*}_{0}\g,\cdots,T^{*}_{0}\g]_0=&[I\oplus I^*,\g\oplus \g^*,\cdots,\g\oplus \g^*]_0\\                                             =&[I,\g,\cdots,\g]_{\g}+[I^*,\g,\cdots,\g]_0+[I,\g,\cdots,\g,\g^*]_0\\
\subseteq& I\oplus I^*=T^{*}_{0}I,
\end{align*}
since
\begin{align*}
[I^*,\g,\cdots,\g]_0(J)=&I^*([J,\g,\cdots,\g]_{\g})\subseteq I^*(J)=0
\intertext{and}
[I,\g,\cdots,\g,\g^*]_0(J)=&\g^*([I,J,\g,\cdots,\g]_{\g})=\g^*(0)=0.
\end{align*}
Then $T^{*}_{0}I$ is a (graded) ideal of $T^{*}_{0}\g$ and so is $T^{*}_{0}J$ in the same way. Hence $T^{*}_{0}\g$ can be decomposed into the direct sum $T^{*}_{0}I\oplus T^{*}_{0}J$ of two nonzero (graded) ideals of $T^{*}_{0}\g$.
\epf

\blem\label{lemma3.1}
Let $(\g,\langle , \rangle_{\theta})$ be a metric first-class $n$-Lie superalgebra of even dimension $m$ over a field $\K$ and $I$ be an isotropic $m/2$-dimensional (graded) subspace of $\g$. Then $I$ is a (graded) ideal of $\g$ if and only if $I$ is abelian.
\elem
\bpf
Since dim$I$+dim$I^{\bot}=m/2+\dim I^{\bot}=m$ and $I\subseteq I^{\bot}$, we have $I=I^{\bot}$.

If $I$ is a (graded) ideal of $\g$, then
$$\langle\g, [\g,\cdots,\g,I,I]_{\g}\rangle_{\theta}=\langle[\g,\cdots,\g,I]_{\g},I\rangle_{\theta}\subseteq \langle I,I\rangle_{\theta}=0,$$
which implies $[\g,\cdots,\g,I,I]_{\g}\subseteq \g^{\bot}=0$.

Conversely, if $[\g,\cdots,\g,I,I]_{\g}=0$, then
$$\langle I,[I,\g,\cdots,\g]_{\g}\rangle_{\theta}=\langle[\g,\cdots,\g,I,I]_{\g},\g\rangle_{\theta}=0.$$ Hence $[I,\g,\cdots,\g]_{\g}\subseteq I^{\bot}=I$. This implies that $I$ is a (graded) ideal of $\g$.
\epf
\bthm\label{theorem3.1}
Let $(\g,\langle , \rangle_{\g})$ be a metric first-class $n$-Lie superalgebra of dimension $m$ over a field $\K$ of characteristic not 2. Then $(\g,\langle , \rangle_{\g})$ is isometric to a $T^{*}$-extension $(T_{\theta}^{*}\g_1,\langle , \rangle_{\theta})$ if and only if $m$ is even and $\g$ contains an isotropic graded ideal $I$ of dimension $m/2$. In particular, $\g_1\cong \g/I$.
\ethm
\bpf
($\Longrightarrow$) Since dim$\g_1$ = dim$\g_1^{*}$, dim$\g$ = dim$T^{*}_{\theta}\g_1=m$ is even. Moreover, it is clear that $\g_1^{*}$ is a graded ideal of  dimension $m/2$ and by the definition of $\langle , \rangle_{\theta}$, we have $\langle \g_1^*,\g_1^* \rangle_{\theta}=0$, i.e., $\g_1^*$ is isotropic.

($\Longleftarrow$) Suppose that $I$ is an $m/2$-dimensional isotropic graded ideal of $\g$. By Lemma \ref{lemma3.1},  $I$ is abelian. Let $\g_1=\g/I$ and $\pi: \g \rightarrow \g_1$ be the canonical projection. Since $\ch \K\neq2$,  we can choose a complement graded subspace $\g_{0}\subseteq\g$ such that $\g=\g_0\dotplus I$ and $\g_0\subseteq \g_0^{\bot}$. Then $\g_0^{\bot}=\g_0$ since dim$\g_0=m/2$.

Denote by $p_{0}$ (resp. $p_1$) the projection $\g \rightarrow \g_0$ (resp. $\g\rightarrow I$) and let $f^*_1$ denote the homogeneous linear map $I \rightarrow \g_1^{*}: z \mapsto f^*_1(z)$, where $f^*_1(z)(\pi(x)):= \langle z,x\rangle_{\g}, \forall x\in \g, \forall z\in I$.

If $\pi(x)=\pi(y)$, then $x-y\in I$, hence $\langle z,x-y\rangle_{\g}\in \langle z,I\rangle_{\g}=0$ and so $\langle z,x\rangle_{\g}=\langle z,y\rangle_{\g}$, which implies $f^*_1$ is well-defined.  Moreover, $f^*_1$ is bijective and $|f^*_1(z)|=|z|$ for all $z\in I$.

In addition, $f^*_1$ has the following property:
\beq\label{eq:f1*&ad*}\begin{split}
&f^*_1([x_1,\cdots,z_k,\cdots,x_n]_{\g})(\pi(y))\\
=&(-1)^{n-k}(-1)^{|z_k|(|x_{k+1}|+\cdots+|x_n|)}\ad^*(\pi(x_1),\cdots,\widehat{\pi(x_k)},\cdots,\pi(x_n))\cdot f^*_1(z_k)(\pi(y)),
\end{split}\eeq
where $x_1,\cdots,x_{k-1},x_{k+1},\cdots,x_n\in \g$, $z_k\in I$.

Define a homogeneous $n$-linear map
\begin{eqnarray*}
\theta:~~~~~ \g_1\times\cdots\times \g_1~~~~&\longrightarrow&\g_1^{*}\\
(\pi(x_1),\cdots,\pi(x_n))&\longmapsto&f^*_1(p_1([x_1,\cdots,x_n]_{\g})),
\end{eqnarray*}
where $x_1,\cdots,x_n\in \g_0.$ Then $\theta$ is well-defined since $\pi|_{\g_0}:\g_0\rightarrow \g_0/I\cong\g/I=\g_1$ is a linear isomorphism and $\theta\in C^1(\g_1,\g_1^*)_{\bar{0}}$.

Now, define the bracket on $\g_1\oplus \g_1^{*}$ by (\ref{eq:bracketofTextension}), then $\g_1\oplus \g_1^{*}$ is an $n$-superalgebra. Let $\varphi$ be a linear map $\g \rightarrow \g_1\oplus \g_1^{*}$ defined by $\varphi(x+z)=\pi(x)+f^*_1(z), \forall x+z\in \g=\g_0\dotplus I. $
 Since $\pi|_{\g_0}$ and $f^*_1$ are linear isomorphisms, $\varphi$ is also a linear isomorphism. Note that
\begin{align*}
&\varphi([x_1+z_1,\cdots,x_n+z_n]_{\g})=\varphi\left([x_1,\cdots,x_n]_{\g}+\sum_{k=1}^n[x_1,\cdots,z_k,\cdots,x_n]_{\g}\right)\\
=&\varphi\left(p_{0}([x_1,\cdots,x_n]_{\g})+p_1([x_1,\cdots,x_n]_{\g})+\sum_{k=1}^n[x_1,\cdots,z_k,\cdots,x_n]_{\g}\right)\\
=&\pi([x_1,\cdots,x_n]_{\g})+f^*_1\left(p_1([x_1,\cdots,x_n]_{\g})+\sum_{k=1}^n[x_1,\cdots,z_k,\cdots,x_n]_{\g}\right)\\
=&[\pi(x_1),\cdots,\pi(x_n)]_{\g_1}+\theta(\pi(x_1),\cdots,\pi(x_n))\\
&+\sum_{k=1}^n(-1)^{n-k}(-1)^{|z_k|(|x_{k+1}|+\cdots+|x_n|)}\ad^*(\pi(x_1),\cdots,\widehat{\pi(x_k)},\cdots,\pi(x_n))\cdot f^*_1(z_k)\\
=&[\pi(x_1)+f^*_1(z_1),\cdots,\pi(x_n)+f^*_1(z_n)]_{\theta}\\
=&[\varphi(x_1+z_1),\cdots,\varphi(x_n+z_n)]_{\theta},
\end{align*}
where we use the definitions of $\varphi$ and $\theta$ and (\ref{eq:f1*&ad*}).
Then $\varphi$ is an isomorphism of $n$-superalgebras, and so $\g_1\oplus \g_1^{*}$ is a first-class $n$-Lie superalgebra.
Furthermore, we have
\begin{align*}
\langle\varphi(x_{0}+z),\varphi(x_{0}'+z')\rangle_{\theta}
&=\langle\pi(x_{0})+f^*_1(z),\pi(x_{0}')+f^*_1(z')\rangle_{\theta}\\
&=f^*_1(z)(\pi(x_{0}'))+(-1)^{|x_0||x_0'|}f^*_1(z')(\pi(x_{0}))\\
&=\langle z,x_{0}'\rangle_{\g}+(-1)^{|x_0||x_0'|}\langle z',x_{0}\rangle_{\g}=\langle x_{0}+z,x_{0}'+z'\rangle_{\g},
\end{align*}
then $\varphi$ is isometric. The relation
\begin{align*}
&\langle[\varphi(x_1+z_1),\cdots,\varphi(x_n+z_n)]_{\theta},\varphi(x_{n+1}+z_{n+1})\rangle_{\theta}\\
=&\langle\varphi([x_1+z_1,\cdots,x_n+z_n]_{\g}),\varphi(x_{n+1}+z_{n+1})\rangle_{\theta}\\
=&\langle[x_1+z_1,\cdots,x_n+z_n]_{\g},x_{n+1}+z_{n+1}\rangle_{\g}\\
=&-(-1)^{(|x_1|+\cdots+|x_{n-1}|)|x_n|}\langle x_n+z_n, [x_1+z_1,\cdots,x_{n-1}+z_{n-1},x_{n+1}+z_{n+1}]_{\g}\rangle_{\g}\\
=&-(-1)^{(|x_1|+\cdots+|x_{n-1}|)|x_n|}\langle \varphi(x_n+z_n), [\varphi(x_1+z_1),\cdots,\varphi(x_{n-1}+z_{n-1}),\varphi(x_{n+1}+z_{n+1})]_{\theta}\rangle_{\theta}
\end{align*}
implies that $(\g_1\oplus \g_1^{*}, \langle ,\rangle_{\theta})$ is a metric first-class $n$-Lie superalgebra.
In this way, we get a $T^*$-extension $T^{*}_{\theta}\g_1$ of $\g_1$ and consequently, $(\g,\langle ,\rangle_{\g})$ and $(T^{*}_{\theta}\g_1,\langle ,\rangle_{\theta})$ are isometric as required.
\epf

Suppose that $\g$ is a first-class $n$-Lie superalgebra and $\theta_{1}$, $\theta_{2}\in Z^1(\g, \g^{*})_{\bar{0}}$ satisfies (\ref{eq:supercyclic}). $T^{*}_{\theta_{1}}\g$ and $T^{*}_{\theta_{2}}\g$ are said to be \textit{equivalent} if there exists an isomorphism of first-class $n$-Lie superalgebras $\phi: T^{*}_{\theta_1}\g\rightarrow  T^{*}_{\theta_2}\g$ such that $\phi|_{\g^*}=\id_{\g^*}$ and the induced map $\bar{\phi}:T^{*}_{\theta_1}\g/\g^{*}\rightarrow T^{*}_{\theta_2}\g/\g^{*}$ is the identity, i.e., $\phi(x)-x\in\g^*$. Moreover, if $\phi$ is also an isometry, then $T^{*}_{\theta_1}\g$ and $T^{*}_{\theta_2}\g$ are said to be \textit{isometrically equivalent}.
\bprop
Suppose that $\g$ is a first-class $n$-Lie superalgebra over a field $\K$ of characteristic not 2 and $\theta_{1}$, $\theta_{2}\in Z^1(\g, \g^{*})_{\bar{0}}$ satisfies (\ref{eq:supercyclic}). Then we have
\begin{enumerate}[(1)]
   \item  $T^{*}_{\theta_{1}}\g$ is equivalent to  $T^{*}_{\theta_{2}}\g$ if and only if $\theta_1-\theta_2\in\delta C^0(\g, \g^{*})_{\bar{0}}$. Moreover, if $\theta_1-\theta_2=\delta\theta'$, then
       \beq\label{eq:induced bilinear form} \langle x, y\rangle_{\theta'}:=\frac{1}{2}\left(\theta'(x)(y)+(-1)^{|x||y|}\theta'(y)(x)\right)\eeq
       becomes a supersymmetric invariant bilinear form on $\g$.
   \item  $T_{\theta_{1}}^{*}\g$ is isometrically equivalent to $T_{\theta_{2}}^{*}\g$ if and only if there is $\theta'\in C^0(\g, \g^{*})_{\bar{0}}$ such that $\theta_1-\theta_2=\delta\theta'$ and the bilinear form induced by $\theta'$ in (\ref{eq:induced bilinear form}) vanishes.
\end{enumerate}
\eprop
\bpf
(1) Let $\phi: T_{\theta_{1}}^{*}\g\rightarrow T_{\theta_{2}}^{*}\g$ be an isomorphism of first-class $n$-Lie superalgebras satisfying $\phi|_{\g^*}=\id_{\g^*}$ and $\phi(x)-x\in \g^*, \forall x\in \g$.
Set $\theta'(x)=\phi(x)-x$. Then $\theta'\in C^0(\g, \g^{*})_{\bar{0}}$ and
\begin{align}
0=&\phi([x_1+f_1,\cdots,x_n+f_n]_{\theta_1})-[\phi(x_1+f_1),\cdots,\phi(x_n+f_n)]_{\theta_2}\notag\\
=&\phi([x_1,\cdots,x_n]_{\g})+\theta_1(x_1,\cdots,x_n)-[x_1+\theta'(x_1)+f_1,\cdots,x_n+\theta'(x_n)+f_n]_{\theta_2}\notag\\
&+\sum_{i=1}^n(-1)^{n-i}(-1)^{|x_i|(|x_{i+1}|+\cdots+|x_n|)}\ad^*(x_1,\cdots,\widehat{x_i},\cdots,x_n)\cdot f_i\notag\\
\begin{split}
=&\theta'([x_1,\cdots,x_n]_{\g})+\theta_1(x_1,\cdots,x_n)-\theta_2(x_1,\cdots,x_n)\\
&-\sum_{i=1}^n(-1)^{n-i}(-1)^{|x_i|(|x_{i+1}|+\cdots+|x_n|)}\ad^*(x_1,\cdots,\widehat{x_i},\cdots,x_n)\cdot \theta'(x_i)
\end{split}\label{eq: theta}\\
=&\theta_1(x_1,\cdots,x_n)-\theta_2(x_1,\cdots,x_n)-\delta\theta'(x_1,\cdots,x_n).\notag
\end{align}

For the converse, suppose that $\theta'\in C^0(\g, \g^{*})_{\bar{0}}$ satisfies $\theta_1-\theta_2=\delta\theta'$. Let $\phi: T_{\theta_{1}}^{*}\g\rightarrow T_{\theta_{2}}^{*}\g$ be defined by $\phi(x+f)=x+\theta'(x)+f$. Then $\phi$ is an isomorphism of first-class $n$-Lie superalgebras such that $\phi|_{\g^*}=\id_{\g^*}$ and $\phi(x)-x\in \g^*, \forall x\in \g$, i.e., $T^{*}_{\theta_1}\g$ is equivalent to  $T^{*}_{\theta_2}\g$.

It's clear that $\langle ,\rangle_{\theta'}$ defined by (\ref{eq:induced bilinear form}) is supersymmetric. Note that
\begin{align*}
&\langle \x\cdot y, z\rangle_{\theta'}+(-1)^{|\x||y|}\langle y, \x\cdot z\rangle_{\theta'}\\
=&\frac{1}{2}\left(\theta'(\x\cdot y)(z)+(-1)^{(|\x|+|y|)|z|}\theta'(z)(\x\cdot y)\right)\\
 &+\frac{1}{2}(-1)^{|\x||y|}\left(\theta'(y)(\x\cdot z)+(-1)^{(|\x|+|z|)|y|}\theta'(\x\cdot z)(y)\right)\\
=&\frac{1}{2}\bigg\{\theta_2(\x,y)(z)-\theta_1(\x,y)(z)+\ad^*(\x)\theta'(y)(z)\\
 &\relphantom{=\frac{1}{2}}+\sum_{i=1}^{n-1}(-1)^{n-i}(-1)^{|x_i|(|x_{i+1}|+\cdots+|x_{n-1}|+|y|)}\ad^*(x_1,\cdots,\widehat{x_i},\cdots,x_{n-1},y)\cdot\theta'(x_i)(z)\bigg\}\\
 &-\frac{1}{2}(-1)^{|y||z|}\ad^*(\x)\cdot \theta'(z)(y)-\frac{1}{2}\ad^*(\x)\cdot \theta'(y)(z)\\
 &+\frac{1}{2}(-1)^{|y||z|}\bigg\{\theta_2(\x,z)(y)-\theta_1(\x,z)(y)+\ad^*(\x)\theta'(z)(y)\\
 &\relphantom{=\frac{1}{2}}+\sum_{i=1}^{n-1}(-1)^{n-i}(-1)^{|x_i|(|x_{i+1}|+\cdots+|x_{n-1}|+|z|)}\ad^*(x_1,\cdots,\widehat{x_i},\cdots,x_{n-1},z)\cdot\theta'(x_i)(y)\bigg\}\\
=&0,
\end{align*}
where we make use of (\ref{eq: theta})=0 and $\theta_1,\theta_2$ satisfying (\ref{eq:supercyclic}).
Then $\langle ,\rangle_{\theta'}$ is invariant.

(2) Let the isomorphism $\phi$ be defined as in (1). Then for all $x+f, y+g\in T^{*}_{\theta_{1}}\g$, we have
\begin{align*}
&\langle\phi(x+f),\phi(y+g)\rangle_{\theta_2}=\langle x+\theta'(x)+f,y+\theta'(y)+g\rangle_{\theta_2}\\
=&\theta'(x)(y)+f(y)+(-1)^{|x||y|}\theta'(y)(x)+(-1)^{|x||y|}g(x)\\
=&2\langle x, y\rangle_{\theta'}+\langle x+f,y+g\rangle_{\theta_1}.
\end{align*}
Thus $\phi$ is an isometry if and only if $\langle , \rangle_{\theta'}=0$.
\epf
\blem\label{lemma3.3}
Let $(V,\langle,\rangle_V)$ be a metric $\Z_2$-graded vector space of dimension $m$ over an algebraically closed field $\K$ of characteristic not 2 and $\g\subseteq gl(V)$ be a Lie superalgebra consisting of nilpotent homogeneous endomorphisms of $V$ such that for each $f\in \g$, the map $f^{+}:V\rightarrow V$ defined by $\langle f^{+}(v),v'\rangle_V=(-1)^{|f||v|}\langle v,f(v')\rangle_V$ is contained in $\g$, too.
Suppose that $W$ is an isotropic graded subspace of $V$ which is stable under $\g$, i.e., $f(W)\subseteq W$ for all $f\in \g$, then $W$ is contained in a maximally isotropic graded subspace $W_{max}$  of  $V$ which is also stable under $\g$ and $\dim W_{max}=[m/2]$. If $m$ is even, then $W_{max}=W_{max}^\bot$. If $m$ is odd, then $W_{max}\subset W_{max}^\bot, \dim  W_{max}^\bot-\dim  W_{max}=1$, and $f(W_{max}^\bot)\subseteq W_{max}$ for all $f\in \g$.
\elem
\bpf
The proof is by induction on $m$. The base step $m = 0$ is
obviously true. For the inductive step, we consider the following two cases.

Case 1: $W \neq 0$ or there is a nonzero $\g$-stable vector $v\in V$(that is, $\g(v)\subseteq\K v$) such that $\langle v,v\rangle_V=0$.

Case 2: $W=0$ and every nonzero $\g$-stable vector $v\in V$ satisfies $\langle v,v\rangle_V\neq 0$.

In the first case $\K v$ is a nonzero isotropic $\g$-stable graded subspace, and $W^{\perp}$ is also $\g$-stable since $\langle w, f(w^\perp)\rangle_V=(-1)^{|f||w|}\langle f^+(w), w^\perp\rangle_V=0$. Now, consider the bilinear form $\langle , \rangle_{V'}$ on the factor graded space $V'=W^{\perp}/W$ defined by $\langle x^\perp+W, y^\perp+W\rangle_{V'}:=\langle x^\perp, y^\perp\rangle_V$, then $V'$ is metric. Denote by $\pi$ the canonical projection $W^{\perp}\rightarrow V'$ and define $f':V'\rightarrow V'$ by $f'(\pi(w^\perp))=\pi(f (w^\perp))$, then $f'$ is well-defined since $W$ and $W^\perp$ are $\g$-stable. Let $\g':= \{f'|f\in\g\}$. Then $\g'$ is a Lie superalgebra. For each $f\in \g$ there is a positive integer $k$ such that $f^{k}=0$, which implies that $(f')^{k}=0$. Hence $\g'$ also consists of nilpotent homogeneous endomorphisms of $V'$. Note that $\g'$ satisfies the same conditions of $\g$. In fact, let $x^\perp$ and $y^\perp$ be two arbitrary elements in $W^{\perp}$. Then by the definition of $\langle , \rangle_{V'}$ we have
\begin{align*}
& \langle (f')^{+}(\pi(x^\perp)), \pi(y^\perp)\rangle_{V'}=(-1)^{|f||x^\perp|}\langle \pi(x^\perp),f'(\pi(y^\perp))\rangle_{V'}\\
=&(-1)^{|f||x^\perp|}\langle \pi(x^\perp),\pi(f(y^\perp))\rangle_{V'}
=(-1)^{|f||x^\perp|}\langle x^\perp,f(y^\perp)\rangle_V\\
=&\langle f^{+}(x^\perp),y^\perp\rangle_V=\langle \pi(f^{+}(x^\perp)),\pi(y^\perp)\rangle_{V'}\\
=&\langle (f^{+})'(\pi(x^\perp)),\pi(y^\perp)\rangle_{V'},
\end{align*}
for arbitrary $f\in \g$, which shows that $(f')^{+}=(f^{+})'\in\g'$ for all $f\in \g$.

Since $\dim V'=\dim W^{\perp}-\dim W=\dim V-2\dim W$, we can use the inductive hypothesis to get a maximally isotropic $\g'$-stable subspace $W'_{max}=W_{max}/W$ in $V'$. Clearly, $\dim W'_{max}$ = $[\frac{\dim V'}{2}]$ = $[\frac{n-2\dim W}{2}]$ = $[n/2]-\dim W$.  For all $x^{\perp}, y^{\perp}\in W_{max}$, the relation $\langle x^{\perp},y^{\perp}\rangle_V=\langle \pi(x^{\perp}),\pi(y^{\perp})\rangle_{V'}=0$ implies that $W_{max}$ is isotropic. Note that $\dim W_{max}=\dim W'_{max}+\dim W=[n/2]$, then $W_{max}$ is maximally isotropic. Moreover, for all $f\in \g$ and $w^{\perp}\in W_{max}$, we have $\pi(f(w^{\perp}))=f'(\pi(w^{\perp}))\in W'_{max}$, which implies $f(w^{\perp})\in W_{max}$. It follows that $W_{max}$ is
$\g$-stable.  This proves the first assertion of the lemma in this case.

In the second case, by Engel's Theorem of Lie superalgebras, there is a nonzero $\g$-stable vector $v\in V$ such that $f(v)=0$ for all $f\in \g$. Clearly, $\K v$ is a nondegenerate $\g$-stable graded subspace of $V$, then $V=\K v\dotplus (\K v)^{\perp}$ and $(\K v)^{\perp}$ is also $\g$-stable since $\langle f((kv)^\perp),v\rangle_V=(-1)^{|f||v|}\langle (kv)^\perp,f^+(v)\rangle_V=(-1)^{|f||v|}\langle (kv)^\perp,0\rangle_V=0, \forall f\in \g$. Now, if $(\K v)^{\perp}=0$, then $V=\K v$ and $\g(V)=0$, hence $\g=0$ and so 0 is the maximally isotropic $\g$-stable subspace, then the lemma follows. If $(\K v)^{\perp}\neq0$, then again by Engel's Theorem of Lie superalgebras there is a nonzero $\g$-stable vector $w\in(\K v)^{\perp}\subseteq V$ such that $f(w)=0$ for all $f \in \g$. It follows that $\g$ vanishes on the two-dimensional nondegenerate subspace $\K v\dotplus \K w$ of $V$. Without loss of generality, we can assume that $\langle v,v\rangle_V=1=\langle w,w\rangle_V$.  Set $\alpha=\langle v,w\rangle_V$, then it is easy to check that the nonzero vector $v+(-\alpha +\sqrt{\alpha^{2}-1})w$ is isotropic and $\g$-stable. This contradicts the assumption of Case 2.

Therefore, the existence of a maximally isotropic $\g$-stable graded subspace $W_{max}$ containing $W$ is proved. If $m$ is even, then dim$W_{max}$=dim$W_{max}^\bot=m/2$; if $m$ is odd, then dim$W_{max}^\bot=\frac{m+1}{2}$ and dim$W_{max}=\frac{m-1}{2}$. Since $\g'$ is nilpotent, there exists a nonzero $\pi(w^\bot)\in V'$ such that $\g'(\pi(w^\bot))=0$. Note that dim$V'$=1, which implies $\g'(V')=0$, so $\g(W_{max}^\bot)\subseteq W_{max}$.
\epf
\bthm\label{thm:2}
Let $(\g,\langle ,\rangle_{\g})$ be a nilpotent metric first-class $n$-Lie superalgebra of dimension $m$ over an algebraically closed field $\K$ of characteristic not 2. If $J$ is an isotropic graded ideal of $\g$, then $\g$ contains a maximally graded ideal $I$ of dimension $[m/2]$  containing $J$. Moreover, if $m$ is even, then $\g$ is isometric to some $T^{*}$-extension of $\g/I$. If $m$ is odd, then $I^\bot$ is abelian and $\g$ is isometric to a nondegenerate graded ideal of codimension 1 in some $T^{*}$-extension of $\g/I$.
\ethm
\bpf
Consider $\ad(\g^{\wedge^{n-1}})=\{\ad\x|\x\in\g^{\wedge^{n-1}}\}$. Then $\ad(\g^{\wedge^{n-1}})$ is a Lie superalgebra. For any $\x\in \g^{\wedge^{n-1}}$, $\ad\x$ is nilpotent since $\g$ is nilpotent. Then the following identity
$$\langle -\ad\x(y), z\rangle_{\g}=(-1)^{|\x||y|}\langle y, \ad\x(z)\rangle_{\g}=(-1)^{|\x||y|}\langle y, \ad\x(z)\rangle_{\g}$$
implies $(\ad\x)^+=-\ad\x\in\g$. Note that $J$ is an $\ad(\g^{\wedge^{n-1}})$-stable graded subspace of $\g$ if and only if $J$ is a graded ideal of $\g$. Then $J$ is an isotropic $\ad(\g^{\wedge^{n-1}})$-stable graded subspace, so there is a maximally isotropic $\ad(\g^{\wedge^{n-1}})$-stable graded subspace $I$ of $\g$ containing $J$ and $\dim I=[m/2]$. Moreover, if $m$ is even, then $\g$ is isometric to some $T^{*}$-extension of $\g/I$ by Theorem \ref{theorem3.1}.

If $m$ is odd, then $\dim I^\bot-\dim I=1$ and $\ad(\g^{\wedge^{n-1}})(I^\bot)\subseteq I$ by Lemma \ref{lemma3.3}. Note that
\begin{align*}
Z(I)=&\{x\in\g|[x,I,\g,\cdots,\g]_{\g}=0\}=\{x\in\g|\langle \g, [x,I,\g,\cdots,\g]_{\g}\rangle_{\g}=0\}\\
    =&\{x\in\g|\langle [I,\g,\cdots,\g]_{\g}, x\rangle_{\g}=0\}=[I,\g,\cdots,\g]_{\g}^\bot=\left(\ad(\g^{\wedge^{n-1}})(I)\right)^{\bot},
\end{align*}
which implies that $I^{\bot}\subset \left(\ad(\g^{\wedge^{n-1}})(I^\bot)\right)^{\bot}=Z(I^{\bot})$, hence $I^{\bot}$ is abelian.

Take any nonzero element $\alpha\notin\g$. Then $\K\alpha$ is a 1-dimensional abelian first-class $n$-Lie superalgebra. Define a bilinear map $\langle ,\rangle_{\alpha}: \K\alpha\times\K\alpha\rightarrow\K$ by $\langle \alpha,\alpha\rangle_{\alpha}=1$. Then $\langle ,\rangle_{\alpha}$ is a nondegenerate supersymmetric invariant bilinear form on $\K\alpha$. Let $\g'=\g\dotplus\K\alpha$. Define
\begin{align*}
[x_1+k_1\alpha,\cdots,x_n+k_n\alpha]_{\g'}=&[x_1,\cdots,x_n]_{\g}\\
\intertext{and}
\langle x+k_1\alpha, y+k_2\alpha\rangle_{\g'}=&\langle x, y\rangle_{\g}+\langle k_1\alpha, k_2\alpha\rangle_{\alpha}.
\end{align*}
Then $(\g',\langle , \rangle_{\g'})$ is a nilpotent metric first-class $n$-Lie superalgebra and $\g$ is a nondegenerate graded ideal of codimension 1 of $\g'$. Since $I^\bot$ is not isotropic and $\K$ is algebraically closed there exists $z\in I^{\bot}$ such that $\langle z, z\rangle_{\g}=-1$. Let $\beta=\alpha+z$ and $I'=I\dotplus\K\beta$. Then $I'$ is an $(m+1)/2$-dimensional isotropic graded ideal of $\g'$.

In fact, for all $x+k_1\alpha+k_1z, y+k_2\alpha+k_2z\in I'$,
\begin{align*}
\langle x+k_1\alpha+k_1z, y+k_2\alpha+k_2z\rangle_{\g'}
=&\langle x+k_1z, y+k_2z\rangle_{\g}+\langle k_1\alpha, k_2\alpha\rangle_{\alpha}\\
=&\langle x, y\rangle_{\g}+\langle x, k_2z\rangle_{\g}+\langle k_1z, y\rangle_{\g}+\langle k_1z, k_2z\rangle_{\g}+k_1k_2\\
=&k_1k_2-k_1k_2=0.
\end{align*}
In light of Theorem \ref{theorem3.1}, we conclude that $\g'$ is isometric to some $T^{*}$-extension of $\g'/I'$.

Define $\Phi:\g'\rightarrow \g/I, x+\lambda\alpha\mapsto x-\lambda z+I$.  Then
\begin{align*}
[\Phi(x_1+\lambda_1\alpha),\cdots,\Phi(x_n+\lambda_n\alpha)]_{\g/I}
=&[x_1-\lambda_1 z+I,\cdots,x_n-\lambda_n z+I]_{\g/I}\\
=&[x_1,\cdots,x_n]_{\g}+I=\Phi([x_1,\cdots,x_n]_{\g})\\
=&\Phi([x_1+\lambda_1\alpha,\cdots,x_n+\lambda_n\alpha]_{\g'}),
\end{align*}
where we use the fact that $I^\bot$ is abelian and $\ad(\g^{\wedge^{n-1}})(I^\bot)\subseteq I$. It's clear that $\Phi$ is surjective and Ker$\Phi=I'$, so $\g'/{I'}\cong \g/I$, hence the theorem follows.
\epf

Now we show that there exists an isotropic graded ideal in every finite-dimensional metric first-class $n$-Lie superalgebra and investigate the nilpotent length of $\g/I$.
\bprop\label{prop:2}
Suppose that $(\g,\langle ,\rangle_{\g})$ is a finite-dimensional metric first-class $n$-Lie superalgebra.
\begin{enumerate}[(1)]
   \item  For any graded subspace $V\subseteq \g$, $C(V):=\{x\in\g|[x,\g,\cdots,\g]_{\g}\subseteq V\}=[\g,\cdots,\g,V^{\bot}]_{\g}^{\bot}$.
   \item  $\g^m=C_m(\g)^{\bot}$, where $C_0(\g)=0, C_{i+1}(\g)=C(C_{i}(\g))$.
   \item  If $\g$ is nilpotent of length $k$, then $\g^i\subseteq C_{k-i}(\g)$.
\end{enumerate}
\eprop
\bpf
The relation
$$\langle C(V),[\g,\cdots,\g,V^{\bot}]_{\g}\rangle_{\g}=\langle [\g,\cdots,\g,C(V)]_{\g}, V^{\bot}\rangle_{\g}\subseteq \langle V, V^{\bot}\rangle_{\g}=0$$
shows that $C(V)\subseteq[\g,\cdots,\g,V^{\bot}]_{\g}^{\bot}$. Notice that
$$\langle [\g,\cdots,\g,[\g,\cdots,\g,V^{\bot}]_{\g}^{\bot}]_{\g}, V^{\bot}\rangle_{\g}=\langle [\g,\cdots,\g,V^{\bot}]_{\g}^{\bot}, [\g,\cdots,\g,V^{\bot}]_{\g}\rangle_{\g}=0,$$
which implies $[\g,\cdots,\g,[\g,\cdots,\g,V^{\bot}]_{\g}^{\bot}]_{\g}\subseteq (V^{\bot})^{\bot}=V$, i.e., $[\g,\cdots,\g,V^{\bot}]_{\g}^{\bot}\subseteq C(V)$. Hence (1) follows.

By induction, (2) and (3) can be proved easily.
\epf

\bthm
Every finite-dimensional nilpotent metric first-class $n$-Lie superalgebra $(\g,\langle ,\rangle_{\g})$ over  an algebraically closed field of characteristic not 2 is isometric to (a nondegenerate ideal of codimension 1 of) a $T^*$-extension of a nilpotent first-class $n$-Lie superalgebra whose nilpotent length is at most a half of the nilpotent length of $\g$.
\ethm
\bpf
Define $J=\sum\limits_{i=0}^{\infty}\g^{i}\cap C_{i}(\g)$. Since $\g$ is nilpotent, the sum is finite. Proposition \ref{prop:2} (2) says $(\g^{i})^{\bot}=C_{i}(\g)$, then $\g^{i}\cap C_{i}(\g)$ is isotropic for all $i\geq 0$.  Since
$$\g^{i}\supseteq \g^{j}\supseteq \g^{j}\cap C_j(g), ~\text{if} ~i<j,$$
we have
$$(\g^{j}\cap C_{j}(\g))^{\bot}\supseteq (\g^{i})^{\bot}=C_{i}(\g)\supseteq C_{i}(\g)\cap \g^{i}, ~\text{if} ~i<j.$$
It follows that
$$\langle\g^{i}\cap C_{i}(\g), \g^{j}\cap C_{j}(\g)\rangle_{\g}=0, ~~\forall i,j\geq0.$$
Therefore $J$ is an isotropic graded ideal of $\g$. Let $k$ denote the nilpotent length of $\g$. Using Proposition \ref{prop:2} (3) we can conclude that $\g^{[(k+1)/2]}\subseteq C_{[(k+1)/2]}(\g)$. This implies that $\g^{[(k+1)/2]}$ is contained in $J$.   By Theorem \ref{thm:2},  there is a maximally isotropic graded ideal $I$ of $\g$ containing $J\supseteq \g^{[(k+1)/2]}$. It means that $\g/I$ has nilpotent length at most $[(k+1)/2]$, and the theorem follows.
\epf

\bre
 Most results of $T^*$-extensions in \cite{B,BL3,LZ} are contained in this section as special cases.
\ere

\end{document}